\documentclass[10pt]{elsarticle}
\usepackage{sectsty}
\usepackage{oubraces}
\subsubsectionfont{\Large}
\usepackage{amsmath}
\usepackage{amsfonts}
\usepackage{multirow}
\usepackage{float}
\usepackage{array}
\usepackage{stix}
\usepackage[a4paper]{geometry}
\usepackage{hyperref}
\usepackage{float}
\usepackage[caption=false]{subfig}
 \newcommand{\dd}{{\rm d}}
\newcommand{\ddF}{{\mathbf{F}}}
\newcommand{\ddp}{{\mathbf{p}}}
\newcommand{\ddq}{{\mathbf{q}}}
\newcommand{\dde}{\mathbf{e}}
\newcommand{\ddR}{{\rm R}}
\newcommand{\ddC}{{\rm C}}
\newcommand{\ddD}{{\rm D}}
\newcommand{\ddc}{{\rm c}}
\newcommand{\ddA}{{\rm A}}
\newcommand{\ddB}{{\rm B}}
\newcommand{\ddO}{{\rm O}}
\journal{Journal of......................}

                    {\endtrivlist}

\begin{document}

\begin{frontmatter}
	
	\title{Analysis of the application of a high order symplectic method in Shardlow's method for dissipative particle dynamics
	}
	\author[a]{Abdolreza Amiri}
	\ead{amiriabdolreza@ymail.com}
	\cortext[mycorrespondingauthor]{Corresponding author}
	
	\address[a]{Department of Mathematics and Statistics, University of Strathclyde, 26
		Richmond Street, Glasgow G1 1XH, Scotland}	

	\begin{abstract}
This study investigates the efficiency and reliability of the modified Shardlow's (M-Shardlow) method for dissipative particle dynamics (DPD). We show that the M-Shardlow method in which for its construction, the second order velocity Verlet method in the Shardlows method to integrate the Hamiltonian part has been replaced by a symplectic fourth order method, improperly uses some parameters. 
By numerical experiments and computing, some important configurational quantities such as configurational temperature and radial distribution function (RDF), the M-Shardlow's method is compared with the Shardlow and ABOBA methods. These results indicate that the new method obtained in this way, even with the proper parameters is too costly in the sense of the CPU-time that is required per each step which makes it an inefficient DPD integrator. Besides, by a comparison of the radial distribution function 
of this method with Shardlow and ABOBA for large time increments, we can observe no considerable improvement in preserving the structure of the system by this new DPD solver.
	\end{abstract}
	
	
\end{frontmatter}

\section{Introduction}
\label{intro}
Dissipative particle dynamics (DPD) is an efficient method to simulate the fluid dynamics of complex fluids. DPD first was proposed by Hoogerbrugge and Koelman \cite{Hoogerbrugge,Koelman}, then, the formulation of the DPD was implemented by Español and Warren \cite{Warren} in 1995. In recent years DPD has been used in many fields of study. It was applied to simulate phase separation in a binary fluid by designating the particles of two different types with a conservative repulsive force between unlike
particles \cite{Coveney}. Also, Schlijper et al. \cite{Schlijper} used DPD  to dilute solutions of polymers by choosing a few of the
particles to represent sections of the polymer chain and adding springs between them. In \cite{Boek} authors
simulate the shear thinning in particulate suspensions. 

Generally, the equations of motion of the DPD particles for $i=1,2, \cdots, N$ are given by
\begin{align}
\dd\ddq_{i}&=m_{i}^{-1}\ddp_{i}\dd t,\nonumber\\
\dd\ddp_{i}&=\ddF_{i}^{\ddC}(\ddq)\dd t+\ddF_{i}^{\ddD}(\ddq,\ddp)\dd t+\dd\ddF_{i}^{\ddR}(\ddq).\label{DPD}
	\end{align}
In the above equations, $\ddF_{i}^{\ddC}(\ddq)$, $\ddF_{i}^{\ddD}(\ddq,\ddp)$ and $\dd\ddF_{i}^{\ddR}(\ddq)$ are respectively, the total conservative, dissipative, and random forces on particle $i$, and are defined as
\begin{subequations}
	\begin{align}
		\ddF_{i}^{\ddC}(\ddq)&=\sum\limits_{j\neq i}\ddF_{ij}^{\ddC}(r_{ij})=-\nabla_{\ddq_{i}} U(\ddq),\label{10a}\\  
		\ddF_{i}^{\ddD}(\ddq,\ddp)&=-\gamma \sum\limits_{j\neq i}w^{\ddD}(r_{ij})(\dde_{ij}\cdot\mathbf{v}_{ij})\dde_{ij},\label{10b}\\
		\dd\ddF_{i}^{\ddR}(\ddq)&=\sigma \sum\limits_{j\neq i}w^{\ddR}(r_{ij})\dde_{ij}\dd\mathbf{W}_{ij}.\label{10c}
	\end{align}
\end{subequations}
where $\dd\mathbf{W}_{ij}$ are independent increments of the Wiener process and $\ddF_{ij}^{\ddC}(r_{ij})$ is the conservative force, and usually is chosen as \cite{Groot}
\begin{eqnarray}
\ddF_{ij}^{\ddC}(r_{ij})=\left\{ \begin{array}{ll} a_{ij}(1-\frac{r_{ij}}{r_{\ddc}})\dde_{ij}, \hspace{1cm}r_{ij}<r_{\ddc};\\  \hspace{0.4cm}0 , \hspace{2.4cm}r_{ij}\geq r_{\ddc},\end{array} \right. \label{HSS}
\end{eqnarray}
 where $a_{ij}$ is the maximum repulsion strength between particles $i$ and $j$, $r_{ij}=\lvert \ddq_{ij}\lvert=\lvert \ddq_{i}-\ddq_{j}\lvert$ is the distance and $r_{\ddc}$ is the cutoff radius. Also, $\dde_{ij}=\ddq_{ij}/r_{ij}$ is the unit vector pointing from particle $i$ to particle $j$ and $\mathbf{v}_{ij}=\ddp_{i}/m_{i}-\ddp_{j}/m_{j}$ is the relative velocity. The parameters $\gamma$, $\sigma$, $w^{\ddD}$ and $w^{\ddR}$ satisfy in the following conditions
\begin{eqnarray}
\sigma^{2}=2\gamma k_{\ddB}T, \hspace{1cm}  w^{\ddD}(r_{ij})=\left(w^{\ddR}(r_{ij}) \right)^{2}.\label{condotions}
\end{eqnarray}
where $k_{\ddB}$ is the Boltzmann constant and $T$ is the the equilibrium temperature. A choice for $w^{\ddR}$ is 
\begin{eqnarray}
w^{\ddR}(r_{ij})=\left\{ \begin{array}{ll} 1-\frac{r_{ij}}{r_{\ddc}}, \hspace{1cm}r_{ij}<r_{\ddc};\\ \hspace{0.4cm} 0 , \hspace{1.4cm}r_{ij}\geq r_{\ddc}.\end{array} \right.
\end{eqnarray}
By the conditions (\ref{condotions}), it is easy to see that the steady-state solution is the Gibbs canonical ensemble
 \begin{eqnarray}
 \rho_{\beta}(\ddq,\ddp)=\dfrac{1}{Z}\exp^{-H(\ddq,\ddp)/k_{\ddB}T},\label{canonical}
 \end{eqnarray}
with $H(\ddq,\ddp)=(\sum\limits_{i}\frac{\ddp_{i}.\ddp_{i}}{2m_{i}}+U(\ddq))/k_{\ddB}T$, and $Z$ is a suitable normalizing constant. Also, $U(\ddq)$ is the potential energy
 \begin{eqnarray}
U(\ddq)=\sum\limits_{i}\sum\limits_{j>i}\varphi(r_{ij}),
 \end{eqnarray}
 where $\varphi(r_{ij})$ is the soft pair potential energy
 \begin{eqnarray}
 \varphi(r_{ij})=\left\{ \begin{array}{ll}\frac{a_{ij}}{2}r_{\ddc} \left(1-\frac{r_{ij}}{r_{c}}\right)^{2}, \hspace{1cm}r_{ij}<r_{\ddc};\\ \hspace{0.4cm} 0 , \hspace{2.75cm}r_{ij}\geq r_{\ddc}.\end{array} \right.
 \end{eqnarray}
 Because of the pairwise (or symmetric) nature of interactions, DPD conserves both angular and linear momenta. So, DPD is an isotropic Galilean-invarriant thermostat which preserve hydrodynamic. 
Consider that when the DPD (\ref{DPD}) preserves the linear momentum, the canonical ensemble with density \cite{Shang} 
\begin{eqnarray}
\rho_{\beta}(\ddq,\ddp)=\dfrac{1}{Z}\exp^{-H(\ddq,\ddp)/k_{\ddB}T}\times \delta \left[\sum\limits_{i}p_{i,x}-\pi_{x}\right] \delta \left[\sum\limits_{i}p_{i,y}-\pi_{y}\right] \delta \left[\sum\limits_{i}p_{i,z}-\pi_{z}\right],\label{ensembel}
\end{eqnarray}
should be replaced (\ref{canonical}), where $\pi=(\pi_{x},\pi_{y},\pi_{z})$ is the linear momentum vector and $\delta$ is the Dirac delta, whose action on a function $\phi$ is $\delta[\phi]=\phi(0)$. If the angular momentum is  conseved, we need to add some additional constraints to (\ref{ensembel}). Note that the ergodicity of the DPD system has only been demonstrated in one dimentioal case with high particle density by Shardlow and Yan \cite{Shardlow}.

Many efforts have been employed to develop accurate and efficient numerical methods to solve (\ref{DPD}). When the time increment is too large, usually these numerical integrators lead to big computational errors. Generally, these efforts to build a new method are in two directions. First, the improvement of time integration schemes, and second the development of the new thermostat. 
The first time-stepping algorithm to solve (\ref{DPD}) was an Euler type method. Then, Groot and Warren \cite{Groot} proposed a modified version of the velocity-Verlet. From that time there have been proposed many numerical integrators to solve the DPD system. These iterative methods aim to produce a higher-order algorithm (usually higher than one) by reducing the additional computations. Many of these iterative methods have stability and efficiency better than the first order Euler and second order velocity-Verlet methods. The focus of most of these methods was on dealing with the time integration and decrease the error for large time increments. Sharlow's method \cite{Shardlow1} is an example of such methods that is an implicit method based on splitting the vector field into a sum of conservative terms and pairwise fluctuation-dissipation terms which allowed for large time
steps when the simulations of the DPD system result in an enormous dissipation or density. 
Also, some of these new methods are Farago and Grønbech-Jensen scheme \cite{Farago} that they construct a new integration method by modifying the St\"{o}mer-Verlet integration (G-JF scheme). Moreover Shang \cite{Shang} recently have presented a
novel splitting method (ABOBA method) which substantially improved the accuracy and efficiency of
DPD simulations for a wide range of friction coefficients.
Based on improving the thermostats, there exist some effective new thermostats, 
Leimkuhler and Shang developed two new stochastic thermostats, pairwise Nos\'{e}-Hoover Langevin (PNHL) thermostat \cite{Leimkuhler} and pairwise adaptive Langevin (PAdL) thermostat \cite{Leimkuhler1}. These two algorithms focus on reducing the effect of dissipative and random forces on numerical error.

Recently, Yamada et al. \cite{Yamada} have proposed that if we substitute a high order symplectic scheme to achieve better accuracy for the Hamiltonian part in Shardlow's method, then the error may decrease for large increments. Consider that the original Shardlow's scheme uses 
the symplectic second-order velocity Verlet method for 
the Hamiltonian part. Therefore, they replaced the velocity Verlet in this method with a fourth-order symplectic method proposed by Yoshida \cite {Yoshida}. They call this modified method M-Shardlow's method. In this paper, we show that they have improperly used the parameters obtained by Yoshida in the fourth-order method to construct their new method, so the results they obtained in this paper are not reliable. Then, the appropriate parameters are used for this method, and theoretically it is shown that with these parameters the M-Shardlow's method with respect to the invariant measure is at least a three order method. Furthermore, by using the concept of numerical efficiency,
we are looking at whether the process used to build M-Shardlow's method can be used as a suitable procedure for building similar methods. We show that even by using the appropriate values of the parameters, this new method does not have the ability to show better performance than the Sharlow or ABOBA methods. The biggest drawback of using this method is the considerable CPU time that is needed for the evaluation of three conservative forces at each iteration compare to Shardlow and the ABOBA methods that use only one evaluation of function.

This paper is organized as follows: In Section 2, first the Shardlow and the ABOBA methods are described, and in this section, by modifying the parameters used in M-Shardlow's method an improvement of this scheme is proposed. Then, theoretically the order of convergence of this new integrator is obtained. In Section 3, to compare this new method especially with the Shardlow and ABOBA methods, a variety of numerical experiments are performed. Also, the numerical efficiency of different methods is investigated in this section. Finally, our findings are summarized in Section 4.  
\section{Numerical methods for DPD}
\label{sec:1} 
 In this section, in addition to a complete description and analysis of the M-Shardlow's method, two methods that have been developed to solve the DPD equation (\ref{DPD}) are briefly described. These two methods are among the best methods proposed for the numerical solution of the DPD system, and it seems that they should be considered a good criterion for benchmarking the newly proposed methods. The first method is Shardlow's method that is a well-known method in solving the system (\ref{DPD}), especially when the DPD system gives rise to a large dissipation or density \cite{Shardlow1}. The second method is the ABOBA method that recently has been introduced by Shang \cite{Shang} and the numerical results show that this method 
 substantially improve the accuracy of the DPD simulations especially for a wide range of friction coefficients. 

\subsection{The Shardlow's splitting method}
Splitting is used widely and as an important approach for solving differential equations. Shardlow's method is an implicit method based on splitting the vector field into a sum of conservative terms (which is well known as the Hamiltonian system $ H $) and pairwise fluctuation-dissipative terms. 

 To describe the Shardlow's method, first, consider the following decomposition for the DPD system (\ref{DPD})
\begin{eqnarray}
\dd\begin{bmatrix} \ddq_{i} \\ \ddp_{i}  \end{bmatrix}=\underset{\ddA}{ \underbrace{\begin{bmatrix} m_{i}^{-1}\ddp_{i} \\ 0  \end{bmatrix}\dd t}}+\underset{\ddB}{ \underbrace{\begin{bmatrix} 0 \\ \ddF_{i}^{\ddC} \end{bmatrix}\dd t}}+\underset{\ddO}{ \underbrace{\begin{bmatrix} 0 \\ \ddF_{i}^{\ddD}\dd t+\dd\ddF_{i}^{\ddR} \end{bmatrix}}},\label{decompose12}
\end{eqnarray}
where $ H=\ddA+\ddB $ and $ \ddO $ part are respectively  corresponding to the Hamiltonian system and Ornstein-Uhlenbeck (OU) process.
To describe the generators of each pieces of the decompsition (\ref{decompose12}), we use the formal  notations that were used in \cite{Serrano1,Fabritiis,Thalmann}. The generators associated with each of the  Hamiltonian system $ H $, and $ \ddO $ parts can be expressed by the following equations
\begin{subequations}
	\begin{align}
	\mathcal{L}_{H}&=\mathcal{L}_{\ddA}+\mathcal{L}_{\ddB},\hspace{0.5cm} \mathcal{L}_{\ddA}=\sum\limits_{i}\frac{\ddp_{i}}{m_{i}}\cdot\nabla_{\ddq_{i}},\hspace{0.5cm} \mathcal{L}_{\ddB}=\sum\limits_{i}\ddF_{i}^{\ddC}\cdot\nabla_{\ddp_{i}}=\sum\limits_{i}\sum\limits_{j> i}\ddF_{ij}^{\ddC}\cdot (\nabla_{\ddp_{i}}-\nabla_{\ddp_{j}}),
	\label{decomp1}\\
	\mathcal{L}_{\ddO}&=\sum\limits_{i}\sum\limits_{j> i}\mathcal{L}_{\ddO_{ij}},\label{decomp3}
	\end{align}
\end{subequations}
where
	\begin{align}
		\mathcal{L}_{\ddO_{ij}}=\left[-\gamma w^{\ddD}(r_{ij})(\dde_{ij}\cdot\mathbf{v}_{ij})+\frac{\sigma^{2}}{2}\left[w^{\ddR}(r_{ij})\right]^{2}\dde_{ij}\cdot(\nabla_{\ddp_{i}}-\nabla_{\ddp_{j}})\right]\dde_{ij}\cdot(\nabla_{\ddp_{i}}-\nabla_{\ddp_{j}}).\label{pairs}
	\end{align}
	is the generator of the $ (i,j) $th interacting pair.
In above equations, $ \mathcal{L}_{\ddO} $ contains the second order derivatives and takes into account the effect of the dissipative and fluctuation terms, and
 $ \mathcal{L}_{\ddA}$ and $\mathcal{L}_{\ddB} $ are respectively the operators for the positions and conservative force.
In Shardlow's method the Verlet method is used to solve the Hamiltonian system. 
Then, it is followed by an implicit method to solve the SDE equation of the fluctuation and dissipation between each interacting pairs.

In Shardlow's method the Br\"{u}nger, Brooks and Karplus (BBK) method \cite{Shardlow, Brooks2} is used successively for integrating the equation generated by (\ref{pairs}) for each interacting pair. 
In \cite{Shardlow} it has been proved that if we looping over all the pairs, and use the BBK method for integrating the SDE equation that is generated by (\ref{pairs}), then a second order method for the 
solution of the whole $ \ddO $ part can be achieved. 
So, if the effect of the BBK method for the $ (i,j) $th interacting pairs be denoted by $ \hat{\mathcal{L}}_{\ddO_{i,j}} $, then the following operator
\begin{eqnarray}
\exp{(\Delta t\hat{\mathcal{L}}_{\ddO})}=\exp{(\Delta t\hat{\mathcal{L}}_{\ddO_{N-1,N}})}\cdots\exp{(\Delta t\hat{\mathcal{L}}_{\ddO_{1,3}})}\exp{(\Delta t\hat{\mathcal{L}}_{\ddO_{1,2}})},\label{pairs12}
\end{eqnarray}
will yield an order $ \Delta t^{2} $ approximation for the OU process \cite{Shardlow}. In  (\ref{pairs12}), $ \Delta t $ is the time integration step and $\exp{(\Delta t \mathcal{L}_{f})}$ denote the flow map of the vector field $f$. 
By these notations the phase space propagation of Shardlow's method that is denoted by DPD-S1, may be expressed as \cite{Shardlow1} 
\begin{eqnarray*}
\exp{(\Delta t\hat{\mathcal{L}}_{DPD-S1})}=\exp{(\Delta t\hat{\mathcal{L}}_{\ddO})}\exp{(\frac{\Delta t}{2}\mathcal{L}_{\ddB})}\exp{(\Delta t\mathcal{L}_{\ddA})}\exp{(\frac{\Delta t}{2}\mathcal{L}_{\ddB})}.
\end{eqnarray*}
Consider that unlike the BBK method that gives a second order approximation for $\ddO$ part, for $\ddA$ and $\ddB$ part we used an exact scheme. The detailed integration steps of the Shardlow's method read as the following:
 \begin{list}{ }{ }
 		\item \textbf{Step 1:}) For each particle $i$ use the velocity Verlet method, that is the following time integration method to obtain $ (\ddq_{i}^{n+1},\ddp_{i}^{n+2/4}) $
 		\begin{subequations}
 			\begin{align}
 			\ddp_{i}^{n+\frac{1}{4}}&=\ddp_{i}^{n}+\dfrac{\Delta t}{2} \ddF_{i}^{\ddC}(\ddq^{n}),\label{V1}\\   
 			\ddq_{i}^{n+1}&=\ddq_{i}^{n}+\Delta t m_{i}^{-1}\ddp_{i}^{n+\frac{1}{4}},\label{V2}\\
 			\ddp_{i}^{n+\frac{2}{4}}&=\ddp_{i}^{n+\frac{1}{4}}+\dfrac{\Delta t}{2} \ddF_{i}^{\ddC}(\ddq^{n+1}).\label{V3}
 			\end{align}
 		\end{subequations}
 	\item \textbf{Step 2:}) For each interacting pair within cutoff radius ($r_{ij}<r_{c}$), updated the momenta as:
\begin{subequations}
	\begin{align}
\ddp_{i}^{n+\frac{3}{4}}&=\ddp_{i}^{n+\frac{2}{4}}-K_{ij}(\dde_{ij}^{n+1}\cdot\mathbf{v}_{ij}^{n+\frac{2}{4}})\dde_{ij}^{n+1}+\mathbf{J}_{ij},\label{p1}\\
\ddp_{j}^{n+\frac{3}{4}}&=\ddp_{j}^{n+\frac{2}{4}}+K_{ij}(\dde_{ij}^{n+1}\cdot\mathbf{v}_{ij}^{n+\frac{2}{4}})\dde_{ij}^{n+1}-\mathbf{J}_{ij},\label{p2}\\
\ddp_{i}^{n+1}&=\ddp_{i}^{n+\frac{3}{4}}+\mathbf{J}_{ij}-\frac{K_{ij}}{1+2K_{ij}}\left[(\dde_{ij}^{n+1}\cdot\mathbf{v}_{ij}^{n+\frac{3}{4}})\dde_{ij}^{n+1}+2\mathbf{J}_{ij}\right],\label{p3}\\
\ddp_{j}^{n+1}&=\ddp_{j}^{n+\frac{3}{4}}-\mathbf{J}_{ij}+\frac{K_{ij}}{1+2K_{ij}}\left[(\dde_{ij}^{n+1}\cdot\mathbf{v}_{ij}^{n+\frac{3}{4}})\dde_{ij}^{n+1}+2\mathbf{J}_{ij}\right],\label{p4}
\end{align}
\end{subequations}
where $K_{ij}=\gamma  w^{\ddD}(r_{ij}^{n+1})\Delta t/2$ and $\mathbf{J}_{ij}=\sigma w^{\ddR}(r_{ij}^{n+1})\dde_{ij}^{n+1} \sqrt{\Delta t} \ddR_{ij}^{n}/2$.
 \end{list}
In equations (\ref{p1})-(\ref{p4}), $\sqrt{\Delta t}\ddR_{ij}$ has been replaced $\dd\mathbf{W}_{ij}$ in the random force, which $\ddR_{ij}$ is a normally distributed variable with zero mean and unit variance .
Note that, the conservative force $\ddF_{i}^{\ddC}$ is computed only once at each iteration, then the updated force is used for the next iteration. Also, at each iteration, the neighbor Verlet list is updated. The interacting pairs can be easily identified from the list in the next iteration to be used to update the momenta.  
\subsection{The ABOBA method}
The ABOBA method is an effective method to solve the DPD system (\ref{DPD}) that is based on the splitting method and 
substantially improve the accuracy and efficiency of the DPD simulation for a wide range of friction coefficients. 
The process that has been described to build the ABOBA method in this section is the same which firstly was done by Shang \cite{Shang}. In the following we explain it briefly. 

First, let us consider the decomposition (\ref{decompose12}) for the DPD system (\ref{DPD}) which the generators of its compartments are given by (\ref{decomp1}) and (\ref{decomp3}) and the  generator corresponding to the $ (i,j) $th interacting pair in the $ \ddO $ part is the equation (\ref{pairs}).

In ABOBA method to solve the system generated by $ \mathcal{L}_{\ddO} $, the $ (i,j) $th SDE for each interacting pair is solved separately by an exact scheme. Obviously, 
 the $ (i,j) $th terms generated by (\ref{pairs}) are the following SDE equations
\begin{subequations}
	\begin{align}
	\dd\ddp_{i}=&-\gamma w^{\ddD}(r_{ij})(\dde_{ij}\cdot\mathbf{v}_{ij})\dde_{ij}\dd t+\sigma w^{\ddR}(r_{ij})\dde_{ij}\dd\mathbf{W}_{ij},\label{SDE3}\\
	\dd\ddp_{j}=&\gamma w^{\ddD}(r_{ij})(\dde_{ij}\cdot\mathbf{v}_{ij})\dde_{ij}\dd t-\sigma  w^{\ddR}(r_{ij})\dde_{ij}\dd\mathbf{W}_{ij}.\label{SDE4}
	\end{align}
\end{subequations}
Let $ \mathbf{v}_{i}=m^{-1}_{i}\ddp_{i} $ be the velocity of the $ i $th particle, then for each interacting pair, subtracting the equation (\ref{SDE4}) from (\ref{SDE3}) yields the following equation
 \begin{align}
 \frac{m_{i}m_{j}}{m_{i}+m_{j}}\dd\mathbf{v}_{ij}=-\gamma w^{\ddD}(r_{ij})(\dde_{ij}\cdot\mathbf{v}_{ij})\dde_{ij}\dd t+\sigma w^{\ddR}(r_{ij})\dde_{ij}\dd\mathbf{W}_{ij}.\label{254}
 \end{align}
 Multiply the unit vector $\dde_{ij}$ on both sides of (\ref{254}) to get 
\begin{align}
 m_{ij}	v_{ij}= -\gamma w^{\ddD}(r_{ij})v_{ij}\dd t+\sigma w^{\ddR}(r_{ij})\dd\mathbf{W}_{ij},\label{DPD2}
 \end{align}
 where $m_{ij}=\frac{m_{i}m_{j}}{m_{i}+m_{j}}$ and  $v_{ij}=\dde_{ij}\cdot\mathbf{v}_{ij}$. (\ref{DPD2}) is an OU process that can be solved exactly. The exact solution of (\ref{DPD2}) is \cite{Kloeden}
 \begin{align}
 v_{ij}(t)=e^{-\tau t}v_{ij}(0)+\frac{\sigma w^{\ddR}(r_{ij})}{m_{ij}}\sqrt{\frac{1-e^{-2\tau t}}{2 \tau}}\ddR_{ij}.\label{O}
 \end{align}
In equation (\ref{O}), $\tau=\frac{\gamma w^{\ddD}(r_{ij})}{m_{ij}} $ and if $v_{ij}(0)$ be the initial relative velocity, then 
 by (\ref{O}), the velocity increment can be obtained as
 \begin{eqnarray}
 \Delta v_{ij}(t)=v_{ij}(t)-v_{ij}(0)=v_{ij}(0)(e^{-\tau t}-1)+\frac{\sigma w^{\ddR}(r_{ij})}{m_{ij}}\sqrt{\frac{1-e^{-2\tau t}}{2 \tau}}\ddR_{ij}.\label{O1}
 \end{eqnarray}
 So, we can update the momenta for each interacting pair as 
 \begin{subequations}\label{O23}
 	\begin{align}
 	\ddp_{i}&\leftarrow \ddp_{i}+m_{ij}\Delta v_{ij}\dde_{ij},\label{O11}\\  
 	\ddp_{j}&\leftarrow\ddp_{j}-m_{ij}\Delta v_{ij}\dde_{ij},\label{O2}
 	\end{align}
 \end{subequations}
 which is coresponding to define the propagator $\exp(\Delta t \mathcal{L}_{\ddO_{i,j}})$ for each interacting pair. 
 
If the scheme (\ref{254})-(\ref{O23}) is employed to solve the corresponding SDE for each interacting pair, then it is easy to see that looping over all the pairs using this exact method, results in an exact integration method for the whole $\ddO$ part. So the overall propagator can be expressed as
   \begin{align}
    \exp(\Delta t \hat{\mathcal{L}}_{\ddO})=\exp{(\Delta t\mathcal{L}_{\ddO_{N-1,N}})}\cdots\exp{(\Delta t \mathcal{L}_{\ddO_{1,3}})}\exp{(\Delta t\mathcal{L}_{\ddO_{1,2}})}.\label{O33}
   	\end{align}
Now, if we use the method, whose the propagator is 
 expressed by  
\begin{eqnarray}
\exp(\frac{\Delta t}{2} \mathcal{L}_{\ddA})\exp(\frac{\Delta t}{2} \mathcal{L}_{\ddB})\exp(\Delta t\hat{\mathcal{L}}_{\ddO})\exp(\frac{\Delta t}{2} \mathcal{L}_{\ddB})\exp(\frac{\Delta t}{2} \mathcal{L}_{\ddA}),\label{ABOBA}
\end{eqnarray}
then, this method is denoted by ABOBA \cite{Shang}.
So, the details of the steps of the ABOBA method  are as the following  \cite{Shang}:
\begin{list}{ }{ }
	\item \textbf{Step 1:}) For particle $i$
		\begin{align*}
	\ddq_{i}^{n+\frac{1}{2}}&=\ddq_{i}^{n}+\dfrac{\Delta t}{2}m_{i}^{-1}\ddp_{i}^{n},\\  
	\ddp_{i}^{n+\frac{1}{3}}&=\ddp_{i}^{n}+\dfrac{\Delta t}{2} \ddF_{i}^{\ddC}(\ddq^{n+\frac{1}{2}}).
\end{align*}
	\item \textbf{Step 2:}) For each interacting pair within cutoff radius ($r_{ij}<r_{c}$), updated the momenta as:
	\begin{align*}
	\ddp_{i}^{n+\frac{2}{3}}&=\ddp_{i}^{n+\frac{1}{3}}+m_{ij}\Delta v_{ij}\dde_{ij}^{n+\frac{1}{2}},\\  
	\ddp_{j}^{n+\frac{2}{3}}&=\ddp_{j}^{n+\frac{1}{3}}-m_{ij}\Delta v_{ij}\dde_{ij}^{n+\frac{1}{2}},
\end{align*}
where
	\begin{eqnarray*}
	\Delta v_{ij}=\left[\dde_{ij}^{n+\frac{1}{2}}\cdot\mathbf{v}_{ij}^{n+\frac{1}{3}} \right](e^{-\tau \Delta t}-1)+\frac{\sigma w^{\ddR}(r_{ij}^{n+\frac{1}{2}})}{m_{ij}}\sqrt{\frac{1-e^{-2\tau \Delta t}}{2 \tau}}\ddR_{ij}^{n},
	\end{eqnarray*}
		with $\tau =\frac{\gamma w^{\ddD}(r_{ij}^{n+\frac{1}{2}})}{m_{ij}} $. 
	\item \textbf{Step 3:}) For each particle $i$
\begin{align*}
	\ddp_{i}^{n+1}&=\ddp_{i}^{n+\frac{2}{3}}+\dfrac{\Delta t}{2} \ddF_{i}^{\ddC}(\ddq^{n+\frac{1}{2}}),\\  
	\ddq_{i}^{n+1}&=\ddq_{i}^{n+\frac{1}{2}}+\dfrac{\Delta t}{2} m_{i}^{-1}\ddp_{i}^{n+1}.
\end{align*}
\end{list}
Note that, in the ABOBA method at each iteration only one conservative force is evaluated. Also, the Verlet list is updated at each iteration, and in the next iteration, the interacting pairs can easily be identified by the Verlet list.
\subsection{M-Shardlow-1 and M-Shardlow-2 splitting methods}
In this section, we study the original  M-Shardlow's method that was presented by  Yamada et al. \cite{Yamada}. Then, we explain that based on the Yoshida \cite{Yamada} scheme to build high order symplectic method, they used improper parameters to construct the modified Shardlow's method. 
Then, the proper parameters are used for M-Shardlow's algorithm and it is shown that with these parameters, the results are much different from the original M-Shardlow's method. From now, we denote the initial M-Shardlow's method by 
M-Shardlow-1 and also, we show by M-Shardlow-2 the improved method that is obtained by using the new parameters.

First, consider that according to the M-Shardlow-1 method, its propagator can be expressed as
 	\begin{align}
 	\exp{(\Delta t\hat{\mathcal{L}}_{{\rm M-Sh-1}})}=\exp{(\Delta t\hat{\mathcal{L}}_{\ddO})}\exp(\Delta t\mathcal{L}_{\mathcal{S}_{4th}}),\label{MShardlow1}
 \end{align}
where $\hat{\mathcal{L}}_{M-Sh-1} $ is the generator of the M-Shardlow-1 scheme and $\hat{\mathcal{L}}_{\ddO}$ is the same BBK generator that was used in Shardlow's method, i.e. for each interacting pair in the $\ddO $ part (OU process), M-Shardlow-1 uses BBK method. Therefore, the propagation of the $ \ddO $ part can be defined as (\ref{pairs12}), and the BBK method is applied successively for integrating each interacting pair in the $ \ddO $ part.
In (\ref{MShardlow1}), $\mathcal{S}_{4th}$ is a fourth order method that has been used for integrating the Hamiltonian part and $ \exp(\Delta t\mathcal{L}_{\mathcal{S}_{4th}}) $ is the propagator of this fourth order method. This fourth order symplectic operator is obtained by the application of Baker-Campbell-Hausdorff (BCH) formula on the second order operator $\exp(\Delta t\mathcal{L}_{\mathcal{S}_{2nd}})$ ($ \mathcal{S}_{2nd} $ shows the second order scheme), namely $\exp(\Delta t\mathcal{L}_{\mathcal{S}_{4th}})$ can be expressed according to the second order operator $\exp(\Delta t\mathcal{L}_{\mathcal{S}_{2nd}})$ as the following
  	\begin{align}
  	\exp(\Delta t\mathcal{L}_{\mathcal{S}_{4th}}):=\exp(x_{1} \Delta t\mathcal{L}_{\mathcal{S}_{2nd}})\exp(x_{0} \Delta t\mathcal{L}_{\mathcal{S}_{2nd}})\exp(x_{1} \Delta t\mathcal{L}_{\mathcal{S}_{2nd}}),\label{BCH1}
  	\end{align}
   provided that, we determine the coefficients $ x_{0} $ and $ x_{1} $ that the exponential series of the right-hand side operator of (\ref{BCH1}) (in the sense of the Taylor exponential series of the operators that has been defined e.g. in \cite{Leimkuhler2}) up to order $4$, coincide with the BCH formula of the Hamiltonian part. In Yoshida scheme second order operators are the basis of building the higher order symplectic methods (see \cite{Yamada}). Let  $ \mathcal{L}_{\ddA}^{\dagger} $ and $\mathcal{L}_{\ddB}^{\dagger}  $ be respectively, the Fokker-Plank operators of the subsystems $ \ddA $ and $ \ddB $ generated by $ \mathcal{L}_{\ddA} $ and $\mathcal{L}_{\ddB}  $. By the BCH formula (\ref{BCH2}), it is obvious that the product 
    \begin{eqnarray}
    \exp(c\Delta t \mathcal{L}_{\ddA}^{\dagger})\exp(d\Delta t\mathcal{L}_{\ddB}^{\dagger})\exp(c\Delta t \mathcal{L}_{\ddA}^{\dagger}),\label{second} \label{product}
    \end{eqnarray}
    is a second order approximation for the effective operator of the Hamiltonian part (the operator $  \exp(\Delta t ( \mathcal{L}_{\ddA}^{\dagger}+\mathcal{L}_{\ddB}^{\dagger})) $), if $ c=\frac{1}{2} $ and $d=1$. Note that, by the BCH expansion (see the formula (\ref{BCH2}) in section \ref{s1}), it is clear that $c$ and $d$ are determined uniquely. Now, if the second order operator (\ref{second}) with $ c=\frac{1}{2} $ and $ d=1 $ be denoted by $\exp( \Delta t\mathcal{L}_{\mathcal{S}_{2nd}})$, again by using the BCH formula we can express the effective operator of the right hand side of (\ref{BCH1}) as a single expansion as \cite{Yoshida}
    \begin{eqnarray}
    \exp(\Delta t (x_{0}+2x_{1})\alpha_{1}+\Delta t^{3}(x_{0}^{3}+2x_{1}^{3})\alpha_{3}+\Delta t^{5}(x_{0}^{5}+2x_{1}^{5})\alpha_{5}+\cdots),\label{BCH3}
     \end{eqnarray}
    where 
      \begin{eqnarray*}
     \alpha_{1}=\mathcal{L}_{\ddA}^{\dagger}+\mathcal{L}_{\ddB}^{\dagger},\hspace{0.2cm} \alpha_{3}=\frac{1}{12}[\mathcal{L}_{\ddB}^{\dagger},\mathcal{L}_{\ddB}^{\dagger},\mathcal{L}_{\ddA}^{\dagger}]-\frac{1}{24}[\mathcal{L}_{\ddA}^{\dagger},\mathcal{L}_{\ddA}^{\dagger},\mathcal{L}_{\ddB}^{\dagger}] ,\hspace{0.2cm}  \alpha_{5}=\frac{7}{5760}[\mathcal{L}_{\ddA}^{\dagger},\mathcal{L}_{\ddA}^{\dagger},\mathcal{L}_{\ddA}^{\dagger},\mathcal{L}_{\ddA}^{\dagger},\mathcal{L}_{\ddB}^{\dagger}]+\cdots. 
      \end{eqnarray*}
      (\ref{BCH3}) gives a fourth order operator for the effective operator of the Hamiltonian part if $x_{0}=-\frac{2^{1/3}}{2-2^{1/3}}$ and $x_{1}=\frac{1}{2-2^{1/3}}$. It is easy to see that these values for $x_{0}$ and $x_{1}$ are determined uniquely.
     But the authors in \cite{Yoshida} used the parameters $x_{0}=\frac{2^{1/3}}{2-2^{1/3}}$ and $x_{1}=\frac{1}{2-2^{1/3}}$ for (\ref{BCH3}), which do not yield a fourth order scheme for the Hamiltonian part. 
     In general, by the exponential representation, the operator $\exp(\Delta t\mathcal{L}_{\mathcal{S}_{4th}})$ can be expressed as the following
      \begin{align}
    \exp(\Delta t\mathcal{L}_{\mathcal{S}_{4th}})&:=\exp(x_{1}\frac{\Delta t}{2}\mathcal{L}_{\ddB})\exp(x_{1}\Delta t\mathcal{L}_{\ddA})\exp(x_{1}\frac{\Delta t}{2}\mathcal{L}_{\ddB})\exp(x_{0}\frac{\Delta t}{2}\mathcal{L}_{\ddB})\exp(x_{0}\Delta t\mathcal{L}_{\ddA})\\ &\exp(x_{0}\frac{\Delta t}{2}\mathcal{L}_{\ddB}) \exp(x_{1}\frac{\Delta t}{2}\mathcal{L}_{\ddB})\exp(x_{1}\Delta t\mathcal{L}_{\ddA}) \exp(x_{1}\frac{\Delta t}{2}\mathcal{L}_{\ddB}).\nonumber
     \end{align}
    So, if the improved method obtained with the new parameters $x_{0}=-\frac{2^{1/3}}{2-2^{1/3}}$ and $x_{1}=\frac{1}{2-2^{1/3}}$ be denoted by M-Shardlow-2, then the details of this method are as the following:
 \begin{list}{ }{ }
 	\item \textbf{Step 1:}) Time integration based on the operator $ \exp(x_{1} \Delta t\mathcal{L}_{\mathcal{S}_{2nd}}) $ for $i$th particle is made as follows to obtain $(\ddq_{i}^{n+1/4},\ddp_{i}^{n+1/4})$: 
 	 \begin{align*}
 		\ddp_{i}^{n+\frac{1}{8}}&=\ddp_{i}^{n}+\dfrac{x_{1}\Delta t}{2} \ddF_{i}^{\ddC}(\ddq^{n}),\\  
 		\ddq_{i}^{n+\frac{1}{4}}&=\ddq_{i}^{n}+x_{1}\Delta t m_{i}^{-1}\ddp_{i}^{n+\frac{1}{8}},\\
 		\ddp_{i}^{n+\frac{1}{4}}&=\ddp_{i}^{n+\frac{1}{8}}+\dfrac{x_{1}\Delta t}{2} \ddF_{i}^{\ddC}(\ddq^{n+\frac{1}{4}}). 
 		  \end{align*}
 		\item \textbf{Step 2:}) Time integration based on the operator $ \exp(x_{0} \Delta t\mathcal{L}_{\mathcal{S}_{2nd}}) $ for $i$th particle is made as follows, to obtain $(\ddq_{i}^{n+2/4},\ddp_{i}^{n+2/4})$: 
 		\begin{align*}
 		\ddp_{i}^{n+\frac{3}{8}}&=\ddp_{i}^{n+\frac{1}{4}}+\dfrac{x_{0}\Delta t}{2} \ddF_{i}^{\ddC}(\ddq^{n+\frac{1}{4}}), \\
 		\ddq_{i}^{n+\frac{2}{4}}&=\ddq_{i}^{n+\frac{1}{4}}+x_{0}\Delta t m_{i}^{-1}\ddp_{i}^{n+\frac{3}{8}},\\
 		\ddp_{i}^{n+\frac{2}{4}}&=\ddp_{i}^{n+\frac{3}{8}}+\dfrac{x_{0}\Delta t}{2} \ddF_{i}^{\ddC}(\ddq^{n+\frac{2}{4}}).
 	\end{align*}
 			\item \textbf{Step 3:}) Time integration based on the operator $ \exp(x_{1} \Delta t\mathcal{L}_{\mathcal{S}_{2nd}}) $ for $i$th particle is made as follows to obtain $(\ddq_{i}^{n+3/4},\ddp_{i}^{n+3/4})$: 
 				\begin{align*}
 				\ddp_{i}^{n+\frac{5}{8}}&=\ddp_{i}^{n+\frac{2}{4}}+\dfrac{x_{1}\Delta t}{2} \ddF_{i}^{\ddC}(\ddq^{n+\frac{2}{4}}), \\
 				\ddq_{i}^{n+\frac{3}{4}}&=\ddq_{i}^{n+\frac{2}{4}}+x_{1}\Delta t m_{i}^{-1}\ddp_{i}^{n+\frac{5}{8}},\\
 				\ddp_{i}^{n+\frac{3}{4}}&=\ddp_{i}^{n+\frac{5}{8}}+\dfrac{x_{1}\Delta t}{2} \ddF_{i}^{\ddC}(\ddq^{n+\frac{3}{4}}). 
 				 	\end{align*}
 				\item \textbf{Step 4:}) For each interacting pair within cutoff radius ($r_{ij}<r_{\ddc}$), update the momenta by the equations (\ref{p1})-(\ref{p2}) and (\ref{p3})-(\ref{p4}).
 \end{list}     
   Note that, because the position is not progressed in the $ \ddO $ part so, in the fourth step, we set $ 	\ddq_{i}^{n+1}=\ddq_{i}^{n+\frac{3}{4}} $. Also, at each iteration and
   at the end of the third step, the neighbor Verlet list is updated, so the interacting pairs can be easily identified from the list and are used for updating the momenta at the fifth step.

     \subsection{Accuracy of the M-Shardlow-2 splitting method to the equilibrium averages}\label{s1}
    In this section, we study the convergence of the M-Shardlow-2 method with respect to the invariant distribution. For this purpose, we need some preliminaries, which are expressed here. It 
    must be noted that the convergence of the Shardlow and ABOBA methods in detail has been explained by Shang \cite{Shang}, so we 
    will not mention them again. 
     
     The comparison of the effective operator of this splitting method with the effective operator of the DPD system is an analysis method for the accuracy of the ergodic averages with respect to the invariant measure that has been used by many authors \cite{Shang,Abdulle,Abdulle1,Leimkuhler12,Leimkuhler13,Talay}. By this analysis, the Baker-Campell-Hausdorff (BCH) formula is used for the comparison of the Taylor expansion of the product of the exponential of the effective operator of  the splitting method with the exponential of the effective operator of the DPD system.
     
     To explain more about this method, consider the splitting method 
     $ \mathcal{L}= \mathcal{L}_{\alpha}+\mathcal{L}_{\beta}+\cdots+\mathcal{L}_{\nu} $, where $ \mathcal{L}_{\alpha},\mathcal{L}_{\beta},\cdots,\mathcal{L}_{\nu} $ are the subsystems of the operator $ \mathcal{L} $. If we denote the exact Fokker-Plank operator by $ \mathcal{L}^{\dagger} $, then $ \mathcal{L}^{\dagger}= \mathcal{L}_{\alpha}^{\dagger}+\mathcal{L}_{\beta}^{\dagger}+\cdots+\mathcal{L}_{\nu}^{\dagger}  $, and the Fokker-Plank operator of each subsystems are shown by $ \mathcal{L}^{\dagger}_{\alpha}, \mathcal{L}^{\dagger}_{\beta},\cdots,\mathcal{L}^{\dagger}_{\nu} $. The effective operator $\hat{\mathcal{L}}^{\dagger}  $ related to the this splitting method can be defined as 
     \begin{align}
     	\exp(\Delta t \hat{\mathcal{L}}^{\dagger})=\exp(\Delta t \mathcal{L}^{\dagger}_{\alpha})\exp(\Delta t \mathcal{L}^{\dagger}_{\beta})\cdots\exp(\Delta t \mathcal{L}^{\dagger}_{\nu}).\label{45}
     \end{align}
    The approximation of the right hand side of (\ref{45}) was denoted by $ 	\exp(\Delta t \hat{\mathcal{L}}^{\dagger}) $. (\ref{45}) can be considered as an approximation of order $ n $ for the operator $ \exp(\Delta t \mathcal{L}^{\dagger}) $, in other words, we have 
      \begin{align}
      \exp(\Delta t \mathcal{L}^{\dagger})=\exp(\Delta t \mathcal{L}^{\dagger}_{\alpha})\exp(\Delta t \mathcal{L}^{\dagger}_{\beta})\cdots\exp(\Delta t \mathcal{L}^{\dagger}_{\nu})+O(\Delta t^{n}).\label{46}
      \end{align}
    where $ O(\Delta t^{n}) $ is the leading error of the splitting method that is obtained by the comparison of the expansion of the right hand side of (\ref{46}) by the BCH formaula with the Taylor expansion of $ \exp(\Delta t \mathcal{L}^{\dagger}) $. For this purpose, consider the following expansion for $ \hat{\mathcal{L}}^{\dagger}. $ 
    \begin{align}
     	\hat{\mathcal{L}}^{\dagger}=\mathcal{L}^{\dagger}+\Delta t \mathcal{L}^{\dagger}_{1}+\Delta t^{2}\mathcal{L}^{\dagger}_{2}+O(\Delta t^{3}),\label{per2}
     \end{align}
      To determine the leading error term in (\ref{per2}), the operators $ \mathcal{L}^{\dagger}_{1} $, $ \mathcal{L}^{\dagger}_{2} $, $ \cdots $ can be computed by the (BCH) expansion.  The expansion (\ref{per2}) characterize the propagation of the density by the effective operator $ \hat{\mathcal{L}}^{\dagger} $. It is clear that the invariant distribution $ \hat{\rho} $ that satisfies in  $  \hat{\mathcal{L}}^{\dagger} \hat{\rho}=0 $ is different from $ \rho_{\beta} $, and we can attribute for every numerical method the pururbed invariant measure $\hat{\rho}$ as an approximation of the targat invariant measure $\rho_{\beta}$, i.e.
     \begin{align}
     	\hat{\rho}=\rho_{\beta}\left[1+\Delta t f_{1}+\Delta t^{2}f_{2}+\Delta t^{3}f_{3}+O(\Delta t^{4})\right],\label{Correction1}
     \end{align}
      where $ f_{1} $, $ f_{2} $, $ \cdots $ are some correction functions.  If we find $ f_{1} $, $ f_{2} $, $ \cdots $, then the invariant distribution which satisfies $  \hat{\mathcal{L}}^{\dagger} \hat{\rho}=0 $ can be determined. 
      This means that the equation $  \hat{\mathcal{L}}^{\dagger} \hat{\rho}=0 $ can be written as
    \begin{align}
     	\left[\mathcal{L}^{\dagger}+\Delta t \mathcal{L}^{\dagger}_{1}+\Delta t^{2}\mathcal{L}^{\dagger}_{2}+O(\Delta t^{3}) \right]\left(\rho_{\beta}\left[1+\Delta t f_{1}+\Delta t^{2}f_{2}+\Delta t^{3}f_{3}+O(\Delta t^{4})\right]\right)=0.\nonumber
      \end{align}
     The exact Fokker-Plank operator preserves the invariant measure, that is, $ \mathcal{L}^{\dagger} \rho_{\beta}=0 $. Also, equating first-order terms in $\Delta t$ yields
     \begin{align}
     	\mathcal{L}^{\dagger}(\rho_{\beta} f_{1})=-\mathcal{L}^{\dagger}_{1}\rho_{\beta}.\label{Per1}
    \end{align}
     Solving (\ref{Per1}) is very difficult, but the BCH formula can be used to find the operator $ \mathcal{L}^{\dagger}_{1} $  and also other operators $\mathcal{L}^{\dagger}_{2}$, $\mathcal{L}^{\dagger}_{3}$, $ \cdots $ in the expansion (\ref{per2}). 
     
     For linear and noncommutative operators $X$ and $Y$ by the BCH expansion, we have the relation 
     \begin{align}
     	\exp(\Delta t X)\exp(\Delta t Y)=\exp(\Delta t Z_{1}).\label{splitting1}
     \end{align}
     where 
     \begin{align}
     	Z_{1}=X+Y+\frac{\Delta t}{2}\left[X,Y\right]+\frac{\Delta t^{2}}{12}\left(\left[X,X,Y\right]-\left[Y,X,Y\right]\right)+O(\Delta t^{3}). \label{BCH}
    \end{align}
     Here we used the notation of the commutator $ \left[X,Y\right]:=XY-YX $. Higher order order commutators like $ \left[Y,\left[X,Y\right]\right] $ was denoted by $ \left[Y,X,Y\right] $. By repeated the application of the BCH formula, we have \cite{Yoshida,Leimkuhler2}
      \begin{align}
     	\exp(\Delta t X) \exp(\Delta t Y)\exp(\Delta t X)=\exp(\Delta t Z_{2}) \label{splitting2}
     \end{align}
     where
     \begin{align}
     	Z_{2}&=2X+Y+\frac{\Delta t^{2}}{6}\left[Y,Y,X\right]-\frac{\Delta t^{2}}{6}\left[X,X,Y\right]\nonumber\\&
     	+
     	\frac{7\Delta t^{4}}{360}\left[X,X,X,X,Y\right]-\frac{\Delta t^{4}}{360}\left[Y,Y,Y,Y,X\right]+\frac{\Delta t^{4}}{90}\left[X,Y,Y,Y,X\right]\label{BCH2}\\&
     	+\frac{\Delta t^{4}}{45}\left[Y,X,X,X,Y\right]-\frac{\Delta t^{4}}{60}\left[X,X,Y,Y,X\right]+\frac{\Delta t^{4}}{30}\left[Y,Y,X,X,Y\right]\nonumber\\&
     	+\cdots.\nonumber
     \end{align}
     So, for the splitting (\ref{splitting1}) that is a nonsymmetric splitting, we have $ \mathcal{L}^{\dagger}_{1}\propto \left[X,Y\right] $, but for the splitting method (\ref{splitting2}), $ \mathcal{L}^{\dagger}_{1}=0 $ and this splitting method achieve a second order method. 

We can use the invariant measure
to compute expectations (averages) of some functions $ \phi $ defined on the phase
space. Such an expectation can be written as    
  \begin{align*}
  	Av_{\phi}=E(\phi)=\int_{\mathcal{D}} \rho_{\beta}(\mathbf{z})\phi(\mathbf{z}) d\mathbf{z} =\lim\limits_{T\rightarrow \infty} \frac{1}{T}\int_{0}^{T}\phi(\mathbf{z})dt.
  	  \end{align*}
  	where  $\mathbf{z}=(\ddp,\ddq)\in\mathcal{D}\subset  \mathbb{\ddR}^{m} $ (for a system of $ N $ particle, typically, $m= 6N $). Any such function is referred to as an observable function (or, simply,
  	 an observable) of the system.
     For an observable $ \phi(\mathbf{z}) $, If $ 	\hat{\rho} $ be the invariant measure  related to the numerical methods, then by (\ref{Correction1}) 
    \begin{align*}
     	\langle \phi \rangle_{\Delta t}&=\int_{\mathcal{D}} \hat{\rho}(\mathbf{z})\phi(\mathbf{z}) d\mathbf{z}=\int_{\mathcal{D}} \rho_{\beta}(\mathbf{z})\left[1+\Delta t f_{1}+\Delta t^{2}f_{2}+\Delta t^{3}f_{3}+O(\Delta t^{4})\right]\phi(\mathbf{z})d\mathbf{z}\\&=\lim\limits_{T\rightarrow \infty} \frac{1}{T}\int_{0}^{T}\left[1+\Delta t f_{1}+\Delta t^{2}f_{2}+\Delta t^{3}f_{3}+O(\Delta t^{4})\right]\phi(\mathbf{z})dt\\&=\langle \phi \rangle + \Delta t \langle \phi f_{1}\rangle + \Delta t^{2} \langle \phi f_{2}\rangle+\Delta t^{3} \langle \phi f_{3}\rangle+ O(\Delta t^{4}).
      \end{align*}
     The notation $ \langle  \rangle_{\Delta t} $ denoted the average of the observable $ \phi(\mathbf{z}) $ in time step $ \Delta t $. It must be noted that the order of convergence also depends on the observable, for instance in some special cases for the observable $ \phi(\mathbf{z}) $ it is possible that $ f_{1}\neq 0 $, but $ \langle \phi f_{1}\rangle =0$, so the method is at least of second order. In certain conditions, higher order convergence to the invariant measure even can be achieved. 
     
     Now, we prove that  M-Shardlow-2 for any observable $ \phi(\mathbf{z}) $ is at least a third order method. As we mentioned  M-Shardlow-2 method uses the BBK successively for integrating each interacting pair in the $ \ddO $ part which is weakly second order consistent \cite{Shardlow}, in other words, the density for each interacting pair can be propagated by an expansion of the form 
       \begin{align}
       \hat{\mathcal{L}}_{\ddO_{i,j}}^{\dagger}=	\mathcal{L}_{\ddO_{i,j}}^{\dagger}+(\Delta t)^{2}\mathcal{L}_{2, \ddO_{i,j}}^{\dagger}+O(\Delta t^{3}),\label{253}
    \end{align}
      where $ \mathcal{L}_{2, \ddO_{i,j}}^{\dagger}$ contains the second order operators in the error of the BBK Fokker-Plank operator when it is used for integrating the $ (i,j) $th interacting pair. Also, by (\ref{pairs12}) and reversal representation of the product which is called Vertauschungssatz \cite{Hairer} for the $\ddO$ part, we have 
       \begin{align}
       \exp{(\Delta t\hat{\mathcal{L}}^{\dagger}_{\ddO})}=\exp{(\Delta t\hat{\mathcal{L}}^{\dagger}_{\ddO_{1,2}})}\exp{(\Delta t\hat{\mathcal{L}}^{\dagger}_{\ddO_{1,3}})}\cdots\exp{(\Delta t\hat{\mathcal{L}}^{\dagger}_{\ddO_{N-1,N}})}.\nonumber
       \end{align}    
       By the BCH formula, the effective operator corresponding to the $ \ddO $ part can be expressed as
       \begin{align}
       \hat{\mathcal{L}}^{\dagger}_{\ddO}=\mathcal{L}^{\dagger}_{\ddO}+\Delta t\overline{\mathcal{L}}^{\dagger}_{1,\ddO}+\Delta t^{2}\overline{\mathcal{L}}^{\dagger}_{2,\ddO}+\Delta t^{3}\overline{\mathcal{L}}^{\dagger}_{3,\ddO}+O(\Delta t^{4}),\label{BCH22}
        \end{align}
       where
       \begin{align}
       \overline{\mathcal{L}}^{\dagger}_{1,\ddO} &=\frac{1}{2}\left[\mathcal{L}^{\dagger}_{\ddO_{1,2}},\mathcal{L}^{\dagger}_{\ddO_{1,3}}\right]+\frac{1}{2}\left[\mathcal{L}^{\dagger}_{\ddO_{1,2}}+\mathcal{L}^{\dagger}_{\ddO_{1,3}},\mathcal{L}^{\dagger}_{\ddO_{1,4}}\right]+\cdots+\frac{1}{2}\left[\mathcal{L}^{\dagger}_{\ddO_{1,2}}+\mathcal{L}^{\dagger}_{\ddO_{1,3}}+\cdots+\mathcal{L}^{\dagger}_{\ddO_{N-2,N}},\mathcal{L}^{\dagger}_{\ddO_{N-1,N}}\right],\nonumber\\
       \overline{\mathcal{L}}^{\dagger}_{2,\ddO} &=\gamma_{2,1}\left[\mathcal{L}^{\dagger}_{\ddO_{1,2}},\mathcal{L}^{\dagger}_{\ddO_{1,2}},\mathcal{L}^{\dagger}_{\ddO_{1,3}}\right]+\gamma_{2,2}\left[\mathcal{L}^{\dagger}_{\ddO_{1,2}},\mathcal{L}^{\dagger}_{\ddO_{1,3}},\mathcal{L}^{\dagger}_{\ddO_{1,3}}\right]\nonumber\\&+\gamma_{2,3}\left[\mathcal{L}^{\dagger}_{\ddO_{1,2}}+\mathcal{L}^{\dagger}_{\ddO_{1,3}},\mathcal{L}^{\dagger}_{\ddO_{1,3}},\mathcal{L}^{\dagger}_{\ddO_{1,3}}\right]+\gamma_{2,4}\left[\mathcal{L}^{\dagger}_{\ddO_{1,2}}+\mathcal{L}^{\dagger}_{\ddO_{1,3}}+\mathcal{L}^{\dagger}_{\ddO_{1,4}},\mathcal{L}^{\dagger}_{\ddO_{1,3}},\mathcal{L}^{\dagger}_{\ddO_{1,3}}\right]+\cdots,\nonumber\\
       \overline{\mathcal{L}}^{\dagger}_{3,\ddO} &=\gamma_{3,1}\left[\mathcal{L}^{\dagger}_{\ddO_{1,2}},\mathcal{L}^{\dagger}_{\ddO_{1,2}},\mathcal{L}^{\dagger}_{\ddO_{1,2}},\mathcal{L}^{\dagger}_{\ddO_{1,3}}\right]+\gamma_{3,2}\left[\mathcal{L}^{\dagger}_{\ddO_{1,2}},\mathcal{L}^{\dagger}_{\ddO_{1,2}},\mathcal{L}^{\dagger}_{\ddO_{1,3}},\mathcal{L}^{\dagger}_{\ddO_{1,3}}\right]\nonumber\\&+\gamma_{3,3}\left[\mathcal{L}^{\dagger}_{\ddO_{1,2}}+\mathcal{L}^{\dagger}_{\ddO_{1,3}},\mathcal{L}^{\dagger}_{\ddO_{1,3}},\mathcal{L}^{\dagger}_{\ddO_{1,3}},\mathcal{L}^{\dagger}_{\ddO_{1,3}}\right]+\gamma_{3,4}\left[\mathcal{L}^{\dagger}_{\ddO_{1,2}}+\mathcal{L}^{\dagger}_{\ddO_{1,3}}+\mathcal{L}^{\dagger}_{\ddO_{1,4}},\mathcal{L}^{\dagger}_{\ddO_{1,3}},\mathcal{L}^{\dagger}_{\ddO_{1,3}},\mathcal{L}^{\dagger}_{\ddO_{1,3}}\right]+\cdots\nonumber\\&+\frac{1}{2}\left[\mathcal{L}^{\dagger}_{\ddO_{1,2}},\mathcal{L}_{2, \ddO_{1,3}}^{\dagger}\right]+\frac{1}{2}\left[\mathcal{L}_{2, \ddO_{1,2}}^{\dagger},\mathcal{L}^{\dagger}_{\ddO_{1,3}}\right]+\frac{1}{2}\left[\mathcal{L}_{2, \ddO_{1,2}}^{\dagger}+\mathcal{L}_{2, \ddO_{1,3}}^{\dagger},\mathcal{L}_{ \ddO_{1,3}}^{\dagger}\right]\nonumber\\&+\frac{1}{2}\left[\mathcal{L}_{2, \ddO_{1,2}}^{\dagger}+\mathcal{L}_{2, \ddO_{1,3}}^{\dagger},\mathcal{L}^{\dagger}_{\ddO_{1,4}}\right]+\cdots+\frac{1}{2}\left[\mathcal{L}^{\dagger}_{2,\ddO_{1,2}}+\mathcal{L}^{\dagger}_{2,\ddO_{1,3}}+\cdots+\mathcal{L}^{\dagger}_{2,\ddO_{N-2,N}},\mathcal{L}_{ \ddO_{N-1,N}}^{\dagger}\right]\nonumber\\&+\frac{1}{2}\left[\mathcal{L}^{\dagger}_{\ddO_{1,2}}+\mathcal{L}^{\dagger}_{\ddO_{1,3}}+\cdots+\mathcal{L}^{\dagger}_{\ddO_{N-2,N}},\mathcal{L}_{ \ddO_{2,N-1,N}}^{\dagger}\right],\nonumber
       \end{align}
       and the coefficients  $ \gamma_{i,j} $, $ i=2,3,\cdots $, $ j=1,2,\cdots $ are determined by the BCH formula. Note that, to write the expansion of $ \overline{\mathcal{L}}^{\dagger}_{3,\ddO} $, the coefficients of the second order error terms in (\ref{253}) have been used. 
       Also, if the effective operator of the fourth order method $S_{4th}$ be denoted by $\mathcal{L}^{\dagger}_{S_{4th}}  $, 
       then 
       \begin{align*}
       \exp(\Delta t \mathcal{L}^{\dagger}_{S_{4th}})=&\exp(x_{1}\frac{\Delta t}{2}\mathcal{L}_{\ddB}^{\dagger})\exp(x_{1}\Delta t\mathcal{L}_{\ddA}^{\dagger})\exp(x_{1}\frac{\Delta t}{2}\mathcal{L}_{\ddB}^{\dagger})\exp(x_{0}\frac{\Delta t}{2}\mathcal{L}_{\ddB}^{\dagger})\exp(x_{0}\Delta t\mathcal{L}_{\ddA}^{\dagger})\\ &\exp(x_{0}\frac{\Delta t}{2}\mathcal{L}_{\ddB}^{\dagger}) \exp(x_{1}\frac{\Delta t}{2}\mathcal{L}_{\ddB}^{\dagger})\exp(x_{1}\Delta t\mathcal{L}_{\ddA}^{\dagger}) \exp(x_{1}\frac{\Delta t}{2}\mathcal{L}_{\ddB}^{\dagger}).
       \end{align*}
       If we set $x_{0}=-\frac{2^{1/3}}{2-2^{1/3}}$ and $x_{1}=\frac{1}{2-2^{1/3}}$, then by using the formula (\ref{BCH3}), we have
       \begin{align}
       \mathcal{L}^{\dagger}_{S_{4th}}=\mathcal{L}_{\ddA}^{\dagger}+ \mathcal{L}_{\ddB}^{\dagger}+\Delta t^{4}\beta_{0}\alpha_{5}(\dfrac{1}{2} \mathcal{L}_{\ddB}^{\dagger},\mathcal{L}_{\ddA}^{\dagger})+O(\Delta t^6),\label{BCH4}
       \end{align}
       where $\beta_{0}=(-\frac{2^{1/3}}{2-2^{1/3}})^{5}+2(\frac{1}{2-2^{1/3}})^{5}$, and we define $ \alpha_{5}(X,Y) $ to be the following operator (see \cite{Yoshida})
       \begin{align}
       \alpha_{5}(X,Y)&=\frac{7}{360}\left[X,X,X,X,Y\right]-\frac{1}{360}\left[Y,Y,Y,Y,X\right]+\frac{1}{90}\left[X,Y,Y,Y,X\right]\nonumber\\&
       +\frac{1}{45}\left[Y,X,X,X,Y\right]-\frac{1}{60}\left[X,X,Y,Y,X\right]+\frac{1}{30}\left[Y,Y,X,X,Y\right]\label{operator}\\&
       +\cdots.\nonumber
       \end{align}
       Finally, we can obtain the effective operator associated with the overall method. First, let $ \mathcal{L}^{\dagger}_{M-Sh-2} $ be the Fokker-Plank operator of the M-Shardlow-2 method, then we have  
     \begin{align}
       \exp(\Delta t\mathcal{L}^{\dagger}_{M-Sh-2})=   \exp(\Delta t\mathcal{L}^{\dagger}_{S_{4th}})\exp(\Delta t\hat{\mathcal{L}}^{\dagger}_{\ddO}).\nonumber	 
       \end{align}
       So, by (\ref{BCH}) the overall BCH expansion can be writen as
        \begin{align}
       \mathcal{L}^{\dagger}_{M-Sh-2}&= \mathcal{L}_{\ddA}^{\dagger}+ \mathcal{L}_{\ddB}^{\dagger}+ \mathcal{L}^{\dagger}_{\ddO}+\Delta t  \mathcal{L}_{1, M-Sh-2}^{\dagger}+\Delta t^{2}  \mathcal{L}_{2, M-Sh-2}^{\dagger}\nonumber	\\ &+\Delta t^{3}  \mathcal{L}_{3, M-Sh-2}^{\dagger}+\Delta t^{4}  \mathcal{L}_{4, M-Sh-2}^{\dagger}+O(\Delta t^{5}),\nonumber
       \end{align}	
       where 
      \begin{align} 
       \mathcal{L}_{1, M-Sh-2}^{\dagger}&=\frac{1}{2}\left[\mathcal{L}_{\ddA}^{\dagger}+ \mathcal{L}_{\ddB}^{\dagger},\mathcal{L}^{\dagger}_{\ddO}\right]+\overline{\mathcal{L}}^{\dagger}_{1,\ddO},\label{523}\nonumber\\
       \mathcal{L}_{2, M-Sh-2}^{\dagger}&=\dfrac{1}{12}\left[\mathcal{L}_{\ddA}^{\dagger}+ \mathcal{L}_{\ddB}^{\dagger},\mathcal{L}_{\ddA}^{\dagger}+ \mathcal{L}_{\ddB}^{\dagger},\mathcal{L}^{\dagger}_{\ddO}\right]+\dfrac{1}{12}\left[\mathcal{L}_{\ddA}^{\dagger}+ \mathcal{L}_{\ddB}^{\dagger},\mathcal{L}^{\dagger}_{\ddO},\mathcal{L}^{\dagger}_{\ddO}\right]+\frac{1}{2}\left[\mathcal{L}_{\ddA}^{\dagger}+ \mathcal{L}_{\ddB}^{\dagger},\overline{\mathcal{L}}^{\dagger}_{1,\ddO}\right]+\overline{\mathcal{L}}^{\dagger}_{2,\ddO},\nonumber\\
       \mathcal{L}_{3, M-Sh-2}^{\dagger}&=\dfrac{1}{24}\left[\mathcal{L}_{\ddA}^{\dagger}+ \mathcal{L}_{\ddB}^{\dagger},\mathcal{L}^{\dagger}_{\ddO},\mathcal{L}^{\dagger}_{\ddO},\mathcal{L}_{\ddA}^{\dagger}+ \mathcal{L}_{\ddB}^{\dagger}\right]+\dfrac{1}{12}\left[\mathcal{L}_{\ddA}^{\dagger}+ \mathcal{L}_{\ddB}^{\dagger},\mathcal{L}_{\ddA}^{\dagger}+ \mathcal{L}_{\ddB}^{\dagger},\overline{\mathcal{L}}^{\dagger}_{1,\ddO}\right]\nonumber\\&+\dfrac{1}{12}\left[\mathcal{L}_{\ddA}^{\dagger}+ \mathcal{L}_{\ddB}^{\dagger},\overline{\mathcal{L}}^{\dagger}_{1,\ddO},\overline{\mathcal{L}}^{\dagger}_{1,\ddO}\right]+\frac{1}{2}\left[\mathcal{L}_{\ddA}^{\dagger}+ \mathcal{L}_{\ddB}^{\dagger},\overline{\mathcal{L}}^{\dagger}_{2,\ddO}\right]+\overline{\mathcal{L}}^{\dagger}_{3,\ddO}.\nonumber
        \end{align}	
       The canonical ensemble is the equilibrium solution of the Hamiltonian system, i.e. $ (\mathcal{L}_{\ddA}^{\dagger}+ \mathcal{L}_{\ddB}^{\dagger})\rho_{\beta}=0 $, also, $ \mathcal{L}_{\ddO_{i,j}}^{\dagger}\rho_{\beta}=0 $, so $ \overline{\mathcal{L}}^{\dagger}_{1,\ddO}\rho_{\beta}=\overline{\mathcal{L}}^{\dagger}_{2,\ddO}\rho_{\beta}=0 $, and we conclude that $  \mathcal{L}_{1, M-Sh-2}^{\dagger}\rho_{\beta} = \mathcal{L}_{2, M-Sh-2}^{\dagger}\rho_{\beta}=0 $. But as the expansion $ \overline{\mathcal{L}}^{\dagger}_{3,\ddO} $ contains  the coefficients of the second order error terms in (\ref{253}) which for that we have $ \overline{\mathcal{L}}^{\dagger}_{3,\ddO}\rho_{\beta}\neq0 $, so $ \mathcal{L}_{3, M-Sh-2}^{\dagger}\rho_{\beta}\neq 0 $. Therefore, for M-Shardlow-2 scheme the correction functions $ f_{1} $ and $ f_{2} $ in (\ref{Correction1}) are zero, and
       the M-Shardlow-2 is at least a third order method.
     \section{Numerical experiments}
     In this section, a set of numerical experiments are implemented to compare the M-Shardlow-2 scheme with M-Sharlow-1 method. Also, we compare this method with the Shardlow and ABOBA methods to see how using a high order symplectic method for the Hamiltonian part can improve the accuracy. Then we employ the concept of numerical efficiency and show that however M-Shardlow-2 method is of a high order, but evaluating three conservative forces at each iteration for this method (compare to Shardlow and ABOBA methods that only evaluate one conservative force at each iteration) makes it somehow ineffective to be used for the simulations. Furthermore, for these methods \lq \lq radial distribution function\rq \rq (RDF) that is also a configurational quantity is  tested in this section.
     \subsection{Simulation detailes}   
     In our numerical experiments, a set of simulations with standard parameters is implemented. First, we apply the same parameters that were used in \cite{Yamada}. The system considered in this paper is a cubic box with length $L=5$ and periodic boundary conditions for all the directions. The DPD parameters were set to be $ m_{i}=1 $, $ r_{\ddc}=1 $, $ a_{ij}=18.75 $, $ \gamma =4.5 $, $ k_{\ddB}T=1 $. Then, we do some other experiments with $a=25$ and different values of $ \gamma $ to study more convergence behavior of the M-Shardlow-2 method. The value of the particle density $ \rho_{\mathrm{DPD}} $ is set according to $ a=75 k_{\ddB}T/\rho_{\mathrm{DPD}} $. Moreover, a system of $ 500 $ identical particles is simulated. The initial position of the particles is distributed uniformly and independently over the cubic box. We choose the initial momenta to be i.i.d (independent identically distributed) normal random variables with mean zero and variance $ k_{\ddB}T $. Each of these simulations is performed for 1000 reduced time units.  
     Verlet neighbor lists \cite{Verlet} were used whenever possible to reduce the computational cost. 
     
 Besides, the computational efficiency of the various methods is tested. For this purpose we use the definition of the \lq \lq numerical efficiency\rq \rq \cite{Leimkuhler,Leimkuhler1}
  that measures the amount of simulation time accessible per unit of computational work, that is
  \begin{eqnarray}
  \text{numerical efficiency}=\frac{\text{Critical Stepsize}}{\text{CPU Time Per Step}},\label{efficiency}
  \end{eqnarray}
  and then scaled it to Shardlow's method (the Shardlow's method is chosen in this paper as the benchmark method). 
  The critical stepsize in (\ref{efficiency}) is determined as the stepsize corresponds to $ 10\% $ relative error in the computed configurational Temperature. 
   The configurational temperature $T_{\ddc}$ is a quantity that measures the temperature and is calculated using the following equation \cite{Braga,Allen,Rugh,Travis}
   \begin{eqnarray}
   k_{\ddB}T_{\ddc}=\frac{\sum_{i}\langle\parallel \nabla_{i}U\parallel^{2}\rangle }{\sum_{i}\langle \nabla_{i}^{2}U\rangle},
  \end{eqnarray}
   where the angle brackets denote the averages. $ \lVert \nabla_{i}U\lVert $ and $ \lVert \nabla_{i}^{2}U\lVert $ are respectively the gradient and Laplacian of the potential energy $U$ with respect to the position of the $ i $th particle.
  Also, in some cases, we employ the kinetic temperature to study some of the results that have been obtained in \cite{Yamada}. The kinetic temperature, $T_{k}$ is defined as 
   \begin{eqnarray}
   k_{\ddB}T_{k}=\frac{1}{d(N-1)}\sum_{i}\frac{\ddp_{i}.\ddp_{i}}{m_{i}},
   \end{eqnarray}
   where $ d $ is the dimension of the physical space (typically $ d =3 $) and $N$ is the number of the particles involved in the simulation. In many previous studies, \cite{Shardlow1,Serrano,Farago,Tang} the kinetic and configurational temperatures have been used to measure the numerical accuracy in the DPD simulations. But in this study, we prefer to use mostly the configurational temperature that is more reliable \cite{Rugh,Butler}, since it can rapidly and accurately track changes in system temperature even when the system is not in global thermodynamic equilibrium. The kinetic temperature depends only on the momenta, and since the canonical ensemble momentum distribution is always Gaussian as in Langevine dynamics, so we are more interested in configurational quantities. Allen recently argued that the configurational temperature should be added to the list of diagnostic tests applied to DPD simulations \cite{Allen2}. The kinetic temperature is used when it is necessary to compare the obtained results for M-Sharlow-1 method in \cite{Yamada} with those of the M-Shardlow-2 method. 
   
   In addition, another important configrational quantity, that is \lq \lq radial distribution function\rq \rq (RDF) \cite{Allen,Frenkel} is calculated. This quantity usually is denoted by  $ g(r) $. RDF is very important  which gives the distribution of distance between particle pairs observed and charectrized the structure of the DPD system, therefore, it can also be considered to check the accuracy of the simulation.
   
    \subsection{Simulation results}   
\begin{figure}[H]
  	\centering
  	\subfloat[Kinetic temperature]{ \includegraphics[width=0.48\textwidth]{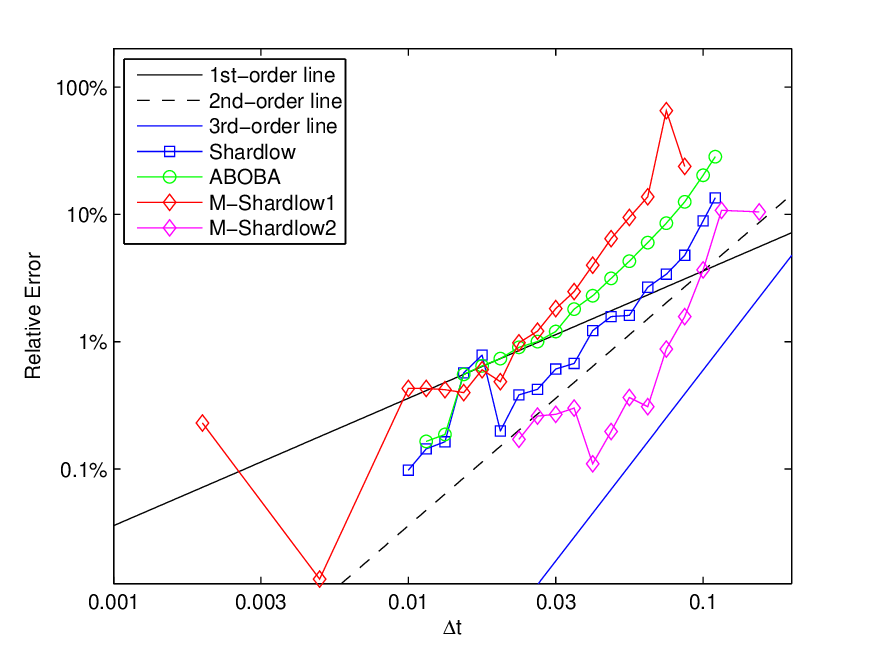}}\hspace{0.2cm}
  	\subfloat[Configurational temperature]{ \includegraphics[width=0.48\textwidth]{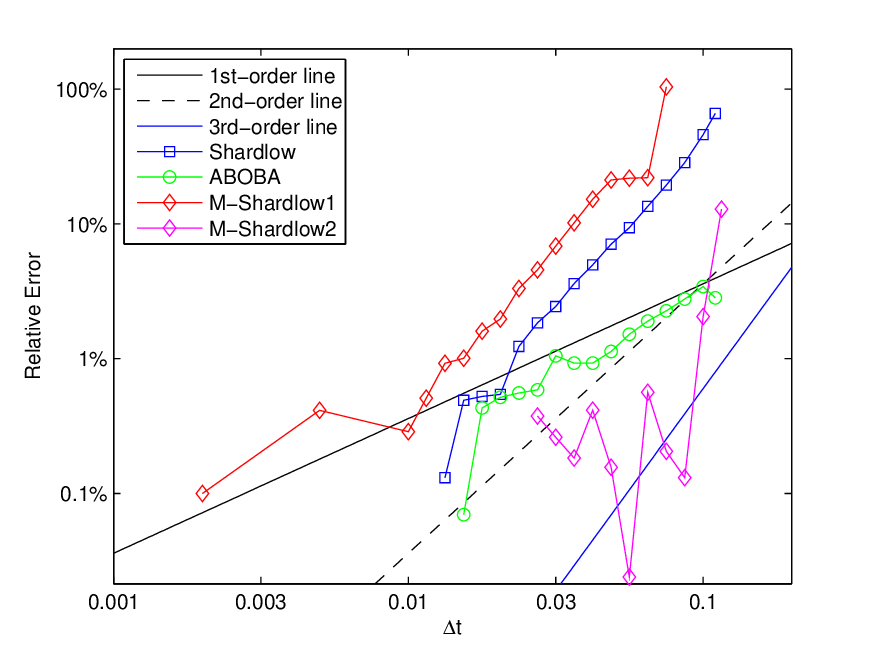}}
  	\caption{The relative error in the computed of the kinetic and configurational temperature for Shardlow, ABOBA, M-Shardlow1, M-Shardlow2 and M-ABOBA methods with $ (a,\gamma)=(18.75,4.5) $. The system was simulated for 1000 reduced time units but only the last $ 80\% $ of the data was collected to compute the static quantity to be sure that the system was well equilibrated. The black and blue solid lines respectively represent the first order and the third order and the black dashed line shows the second order convergence to the invariant measure.} \label{fig1}
  \end{figure}

  The configurational temperature control of different methods is investigated in this section. We compare the configurational temperature of the M-Shardlow-1 method with the M-Shardlow-2 method. Figure \ref{fig1} shows the relative error in the computed of the kinetic and configurational temperature when the simulations are implemented with parameters $ (a_{ij},\gamma)=(18.75,4.5) $. These parameters are the same parameters that have been used by the authors in \cite{Yamada}. The simulations were implemented for the time increments ranging between $0.002$ and $ 0.12 $. The stepsizes began at around $ \Delta t = 0.002 $ and increased
 until the method shows significant relative errors or became unstable. For Shardlow, M-Shardlow-1 and ABOBA methods the error in both kinetic and configurational temperatures tend to increase after $\Delta t\geq 0.01$. In these three methods, the error manifests almost second order convergence at the range between $ 0.02$ and $ 0.1 $. Similar to the results that obtained in \cite{Yamada} 
  for smaller stepsizes ($ \Delta t <0.02 $), it is found that in the M-Shardlow-1 method, the increase in the error is very slow and manifests first order convergence or even slower, however for other methods also slow convergence for small increments can be expected. These results as the theorotical ones show that as the M-Shardlow-1 method is using improper parameters, so the order of convergence higher than two for this method  can not be expected for this method. 
  
  For the M-Shardlow-2 method a very gradual increase in error can be observed with increasing time increments, however when approaching $ \Delta t=0.1 $, it manifest some high-order convergence, but we can not infer three order convergence explicitly for the M-Shardlow-2 method by this simulation. So, we set a various of other experiments with $a=25$ and $\gamma=4.5,40.5,200,450$. It is obvious that by the increase of the parameter $ \gamma $ we expect the error to grow more rapidly for large time increments and is more convenient to study the convergence behavior. 
  
\begin{figure}[H]
				 	\centering
				 	\subfloat[$\gamma=4.5$]{ \includegraphics[width=0.48\textwidth]{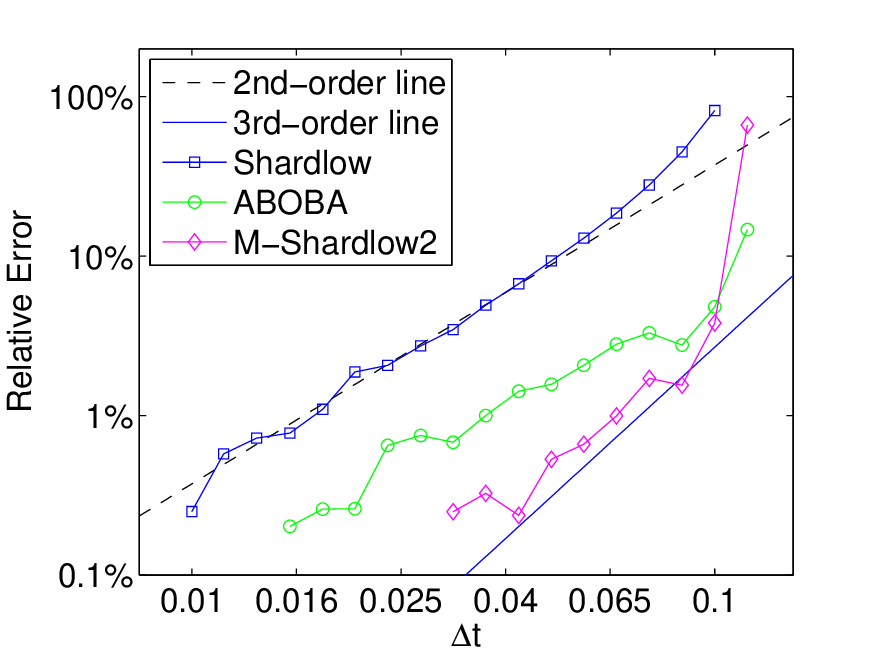}}\hspace{0.2cm}
				 	\subfloat[$\gamma=40.5$]{ \includegraphics[width=0.48\textwidth]{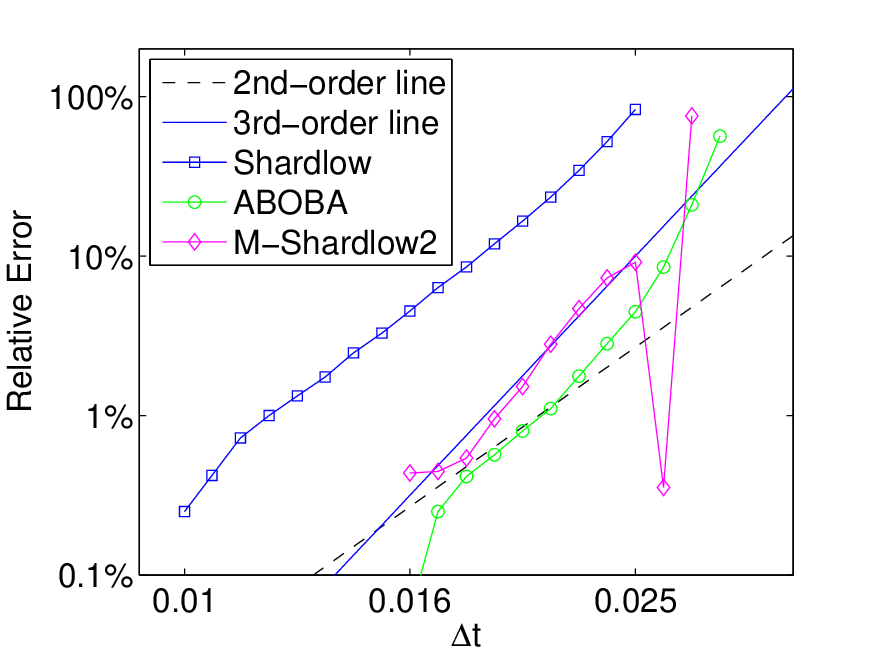}}\\
				 	\subfloat[$\gamma=200$]{ \includegraphics[width=0.48\textwidth]{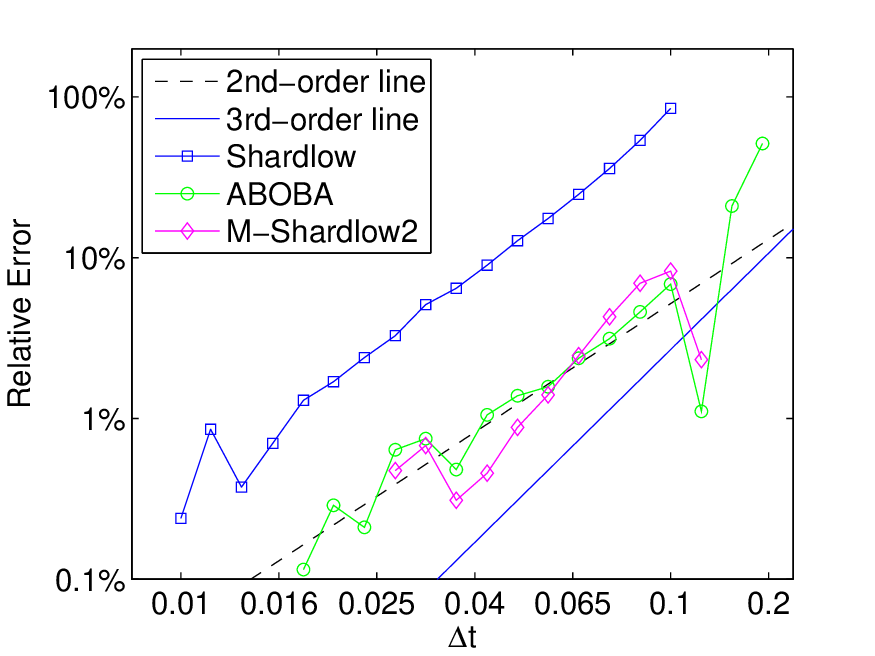}}\hspace{0.2cm}
				 	\subfloat[$\gamma=450$]{ \includegraphics[width=0.48\textwidth]{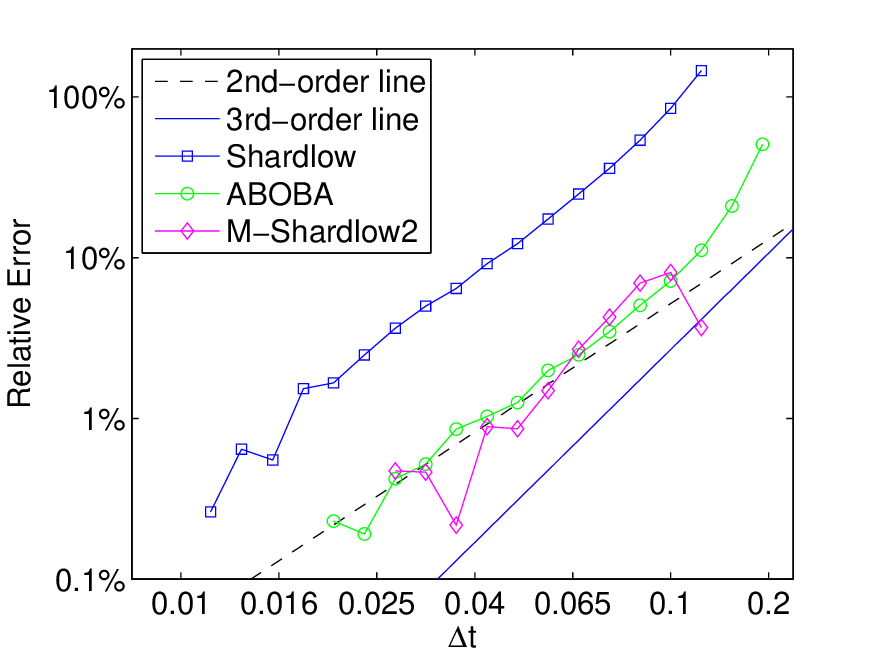}}
				 	\caption{Double logarithmic of the relative error in the computed configurational temperature in Shardlow, ABOBA, M-Shardlow2, and M-ABOBA methods for $\gamma=4.5,40.5,200,450$ ($a=25$). The format of the plots is the same as in figure \ref{fig1}, except that the stepsize starts at $ \Delta t=0.01 $.
				 		}\label{fig2}
				 \end{figure}

 Figure \ref{fig2} shows the computed configurational temperature in Shardlow, ABOBA, and M-Shardlow-2 methods for $ a=25 $ and different values of $\gamma$. It was mentioned earlier in this section that we prefer only to compare the configurational temperature, so we excluded the kinetic temperature in this case. The time increments ranging between $0.01$ and $ 0.15 $. As expected and has been shown recently by Shang \cite{Shang} that both ABOBA and Shardlow's methods successfully for high frictions implement the simulations and show a second order convergence. For M-Shardlow-2 method, the error increases very slowly for $ \Delta t <0.04 $, but after $ \Delta t =0.04 $ a third order convergence can easily be observed. In these figures, we didn't include M-Shardlow-1 method, since it uses improper parameters and usually the results are not appropriate and also are worst than the other three methods.  

	In Tables 1-4, the \lq \lq numerical efficiency\rq \rq, defined as the ratio of the critical stepsize and CPU time per step, of each method for $a=25$ and different values of $ \gamma $ has been measured. In our computations, we consider the critical stepsize when the relative error be $ 10\% $ in the computed configurational temperature according to Figure \ref{fig2}. The CPU time per each step obtained, then scaled to Shardlow's method. To obtain the CPU time for one step, the simulations performed for 10000 successive time steps, then the results averaged over 10000 steps. The CPU time is in milliseconds and the simulations were performed on an Intel Core processor i7-4700HQ with a CPU of 2.40GHz and 8.0 GB of RAM. Physical quantities should not be taken into account when calculating the CPU time. Shardlow's method was chosen as the benchmark method, and other methods have been compared with this method. As mentioned in the previous section, the Verlet list is used in computations whenever possible. For the Shardlow, M-Shardlow-2, and ABOBA methods, when the conservative forces are computed, it is possible to list the interacting pairs and used them to identify the interacting pairs for updating the momenta in the dissipative-fluctuation part. M-Shardlow-2 method in all cases outperforms the M-Shardlow-1 method, so M-Shardlow-1 method was not included in these tables.
			
		In Table 1 the numerical efficiency for $ (a,\gamma)=(25,4.5) $ has been measured for different methods. This table shows that the ABOBA method outperforms the other three methods, while M-Shardlow-2  has achieved the performance of about $ 71.55 \%$. However, when we use a high order symplectic method for the Hamiltonian part it increases the stability threshold considerably, but for the M-Shardlow-2 method, three conservative forces are needed to be evaluated at each iteration, so the computational costs outweigh the considerable improvement in the threshold stability. For this reason, one can see that the M-Shardlow-2 does not improve the numerical efficiency with respect to Shardlow's method, and even the numerical efficiency of the M-Shardlow-2 method has decreased compared to the Shardlow's methods and is significantly lower than the ABOBA method.  
		
\begin{figure}[H]
\centering
\subfloat[Shardlow]{ \includegraphics[width=0.48\textwidth]{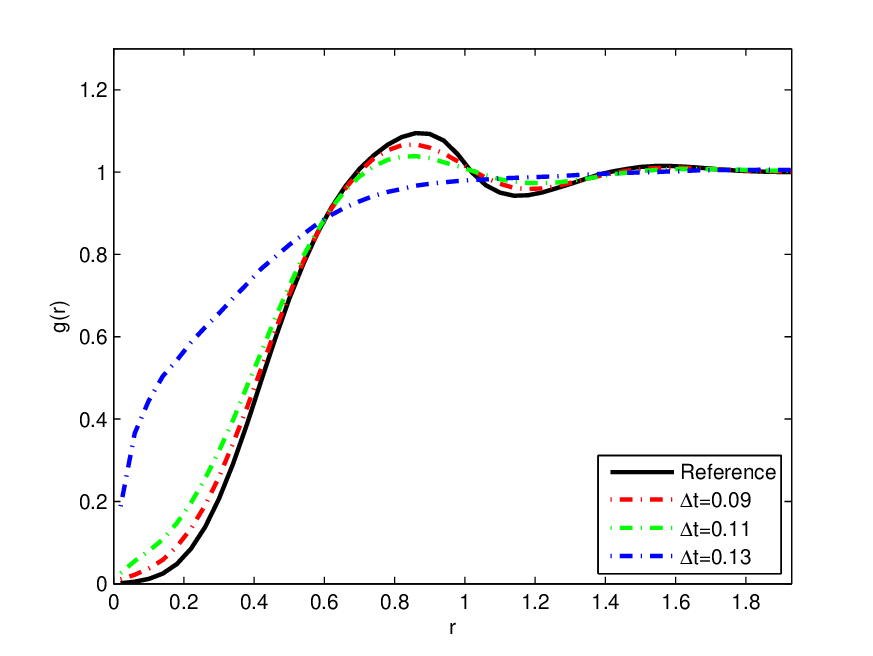}}\hspace{0.2cm}
\subfloat[ABOBA]{ \includegraphics[width=0.48\textwidth]{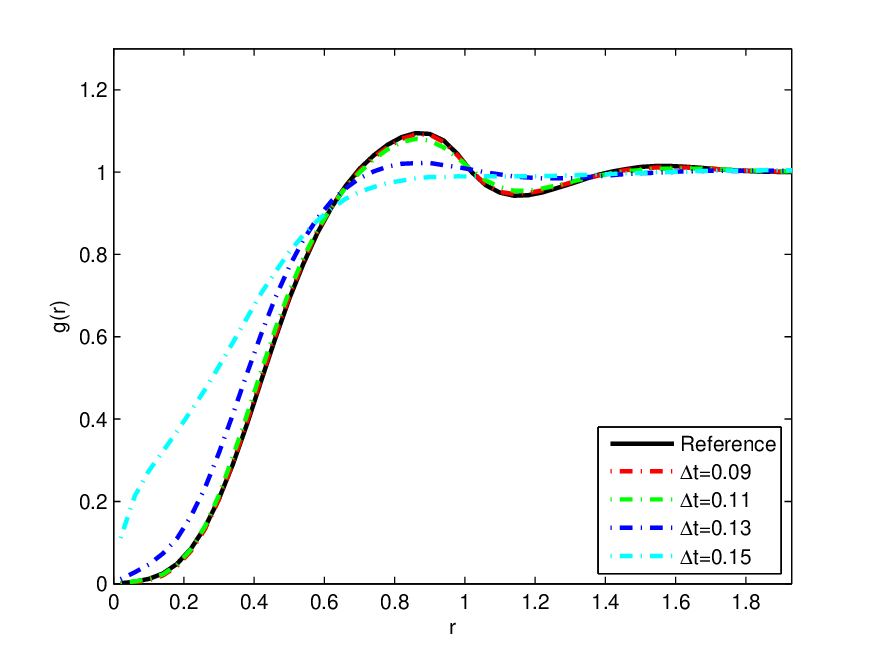}}\\
\subfloat[M-Shardlow-1]{ \includegraphics[width=0.48\textwidth]{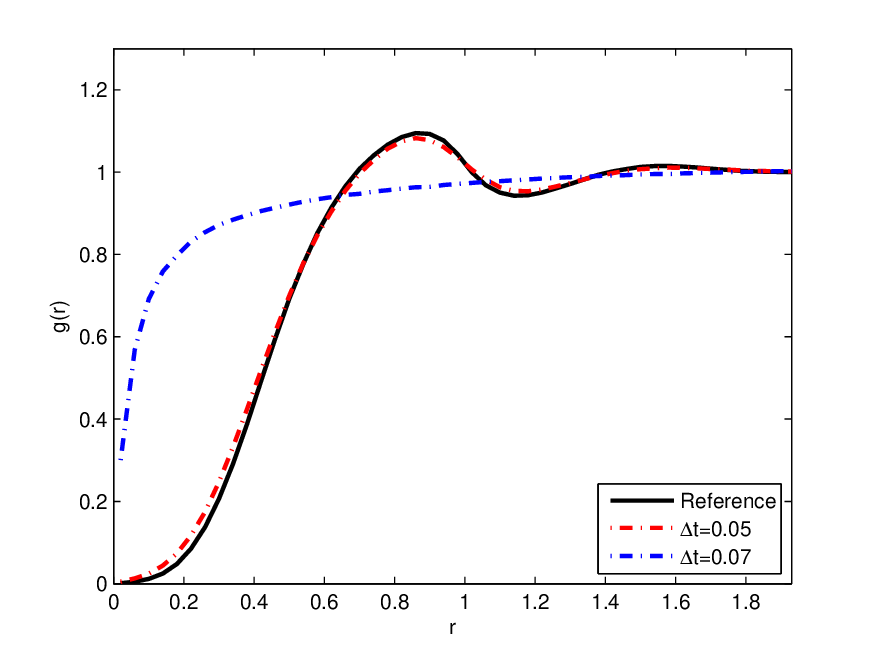}}\hspace{0.2cm}
\subfloat[M-Shardlow-2]{ \includegraphics[width=0.48\textwidth]{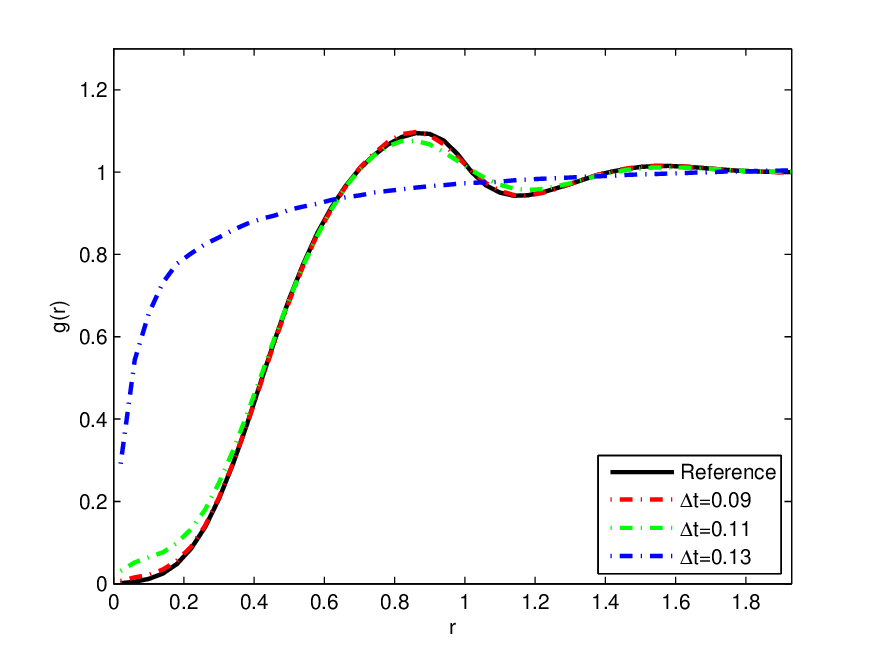}}
\caption{Comparison of the radial distribution function (RDF), $ g(r) $, for different methods with a friction coefficient of $ \gamma=4.5 $ and constant repulsion $ a=25 $. In all the methods the stepsize start at $ \Delta t=0.09 $ and RDFs have been plotted until large artifacts observed. The solid black line is the refrence solution obtained by using the Shardlow's method with a very small stepsize 
of $ \Delta t=0.001 $. The dashed colour lines correspond to different stepsize as indicated.}\label{fig3}
\end{figure}

				  \begin{table}[H]
				  	\begin{center}
				  		\scalebox{1.1}{
				  			\begin{tabular}{lccccccccccccccr}
				  				\hline
				  				Method        &  & & Critical stepsize& CPU time&Scaled efficiency\\
				  				\hline
				  				Shardlow 	&& & 0.051   &13.0&100.0$\%$ \\ 
				  				M-Shardlow-2&& & 0.096   &34.2&  71.55$\%$ \\ 
				  				ABOBA	& & & 0.1132  &11.50& 250.91 $\%$ \\ 
				  				\hline
				  			\end{tabular}}
				  		\end{center}
				  		\vspace{-0.6cm}
				  		\caption{Comparison of the \lq \lq numerical efficiency\rq \rq of the Shardlow, ABOBA and M-Shardlow-2 methods with $(\gamma,a) =(4.5,25)$. The numerical efficiency of each method was scaled to that of Shardlow's method.}
				  	\end{table}
				  	In Tables 2 the numerical efficiency has been computed for $ (a,\gamma)=(25,40.5) $. By this table again we can see a slight improvement in the performance of the M-Shardlow-2
				  	method compare to the Shardlow's method. In this case the usage of the M-Shardlow-2 method is significantly less justified than the ABOBA method, since by this table the efficiency of the ABOBA method is about $ 298\% $ better than the M-Shardlow-2 method.
				  	\begin{table}[H]
				  		\begin{center}
				  			\scalebox{1.1}{
				  				\begin{tabular}{lccccccccccccccr}
				  					\hline
				  					Method        &  & & Critical stepsize& CPU time&Scaled efficiency\\
				  					\hline
				  					Shardlow 	&& & 0.0455   &13.0&100.0$\%$ \\ 
				  					M-Shardlow-2& & & 0.1192   &34.0&  100.84$\%$ \\ 
				  					ABOBA	& & & 0.12  &11.50&  298.14$\%$ \\ 
				  					\hline
				  				\end{tabular}}
				  			\end{center}
				  			\vspace{-0.6cm}
				  			\caption{Comparison of the \lq \lq numerical efficiency\rq \rq of the Shardlow, ABOBA and M-Shardlow-2 methods with $(\gamma,a) =(40.5,25)$. The numerical efficiency of each method was scaled to that of Shardlow's method.}
				  		\end{table}
				  	Similar results to Table 2 have been obtained for these three methods, when the parameters $ (a,\gamma) =(25,200)$ and $ (a,\gamma) =(25,450)$ are used in the simulations. The results for the performance of the various methods for these parameters are reported in Tables 3 and 4. By these tables again, no significant improvement in the performance of 
				  	the M-Shardlow-2 method is observed compared to other methods. 
				  		\begin{table}[H]
				  			\begin{center}
				  				\scalebox{1.1}{
				  					\begin{tabular}{lccccccccccccccr}
				  						\hline
				  						Method        &  & & Critical stepsize& CPU time&Scaled efficiency\\
				  						\hline
				  						Shardlow 	&& & 0.0443   &12.9&100.0$\%$\\ 
				  						M-Shardlow-2& & & 0.118   &34.0&  101.06$\%$ \\ 
				  						ABOBA	& & &0.114  &11.50&  288.66 $\%$ \\ 
				  						\hline
				  					\end{tabular}}
				  				\end{center}
				  				\vspace{-0.6cm}
				  				\caption{Comparison of the \lq \lq numerical efficiency\rq \rq of the Shardlow, ABOBA and M-Shardlow-2 methods with $(\gamma,a) =(200,25)$. The numerical efficiency of each method was scaled to that of Shardlow's method.}	
				  			\end{table}
				  			\begin{table}[H]
				  				\begin{center}
				  					\scalebox{1.1}{
				  						\begin{tabular}{lccccccccccccccr}
				  							\hline
				  							Method        &  & & Critical stepsize& CPU time&Scaled efficiency\\
				  							\hline
				  							Shardlow 	&& & 0.044   &12.8&100.0$\%$\\ 
				  							M-Shardlow-2& & & 0.117   &34.7& 109.59$\%$ \\ 
				  							ABOBA	& & & 0.112  &11.48&277.11$\%$ \\ 
				  							\hline
				  						\end{tabular}}
				  					\end{center}
				  					\vspace{-0.6cm}
				  					\caption{Comparison of the \lq \lq numerical efficiency\rq \rq of the Shardlow, ABOBA  and M-Shardlow-2 methods with $(\gamma,a) =(450,25)$. The numerical efficiency of each method was scaled to that of Shardlow's method.}	
				  				\end{table}
				  		Figure \ref{fig3} compares the RDF of various methods for the friction coefficient $ \gamma=4.5 $ and repulsion $ a=25 $. These results are quite consistent with the results obtained for configurational temperature. In these figures the reference RDF that was indicated by the solid black lines obtained by using Shardlow's method with a very small stepsize of $ \Delta t=0.001 $. Except for the M-Shardlow-1 method, the RDFs are plotted with increments equal and bigger than $ \Delta t=0.09 $, and for each one, the simulations continued until large artifacts be noticed. RDF for the M-Shardlow-1 method is plotted for smaller stepsizes. As it was expected for the Shardlow1's method based on our findings that we obtained for the configurational temperature, it is very unstable after $ 0.065 $ and the RDF shows the same results. The Shardlow and M-Shardlow-2 methods both show artifacts at around $ \Delta t=0.13 $ and the RDFs heavily destroyed at this stepsize. The ABOBA method even shows better performance in conserving the structure of the system and a very small distortion of the reference RDF for stepsize of $ 0.13 $ is observed. This shows that using a fourth order method for Shardlow's method didn't help too much to preserve the structure of the DPD system. Even for this friction coefficient, the ABOBA method shows better control of the structure of the system than the M-Shardlow-2 method.
				  			
\begin{figure}[H]
\centering
		\subfloat[Shardlow]{ \includegraphics[width=0.48\textwidth]{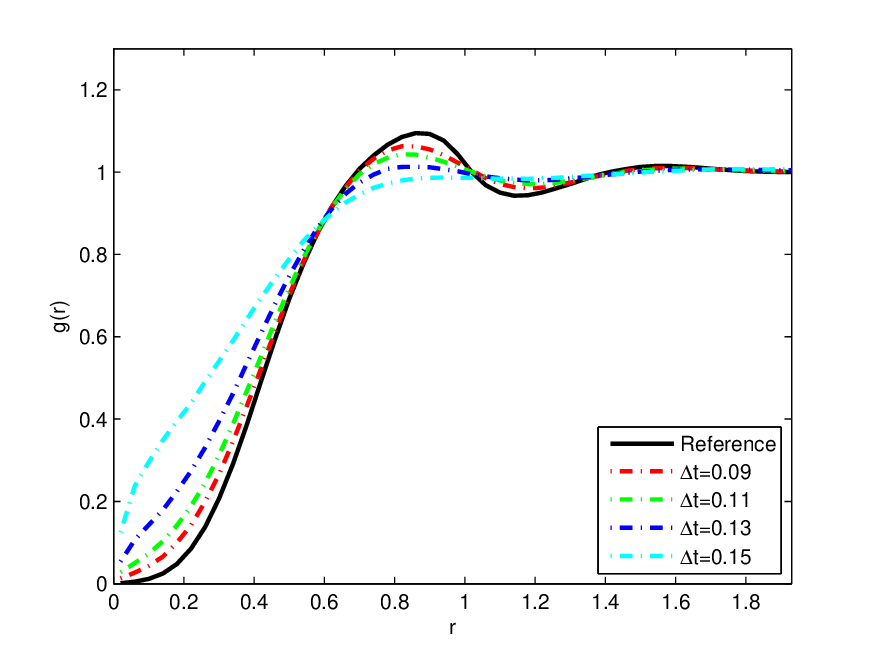}}\hspace{0.2cm}
	\subfloat[ABOBA]{ \includegraphics[width=0.48\textwidth]{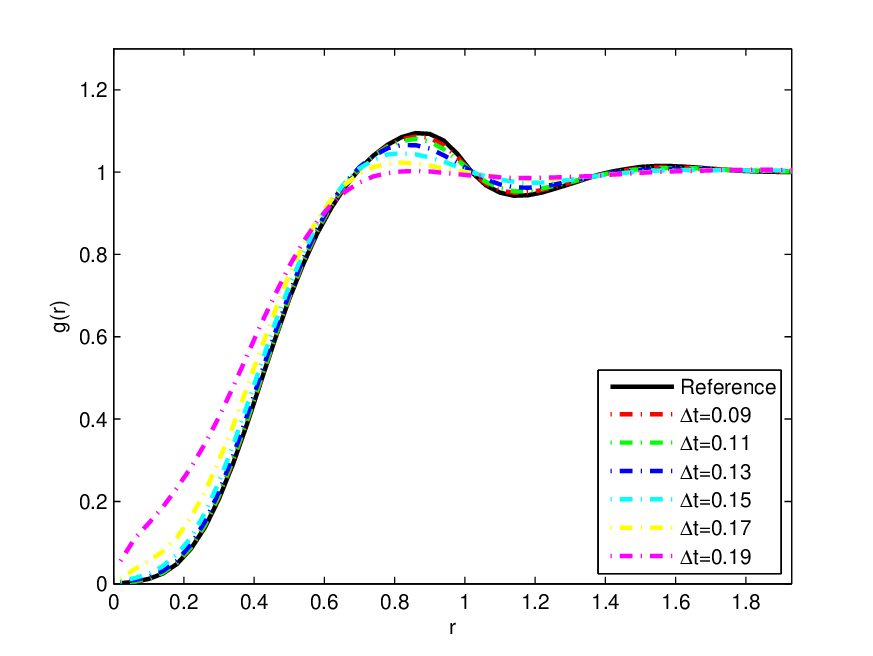}}\\
	\subfloat[M-Shardlow-1]{ \includegraphics[width=0.48\textwidth]{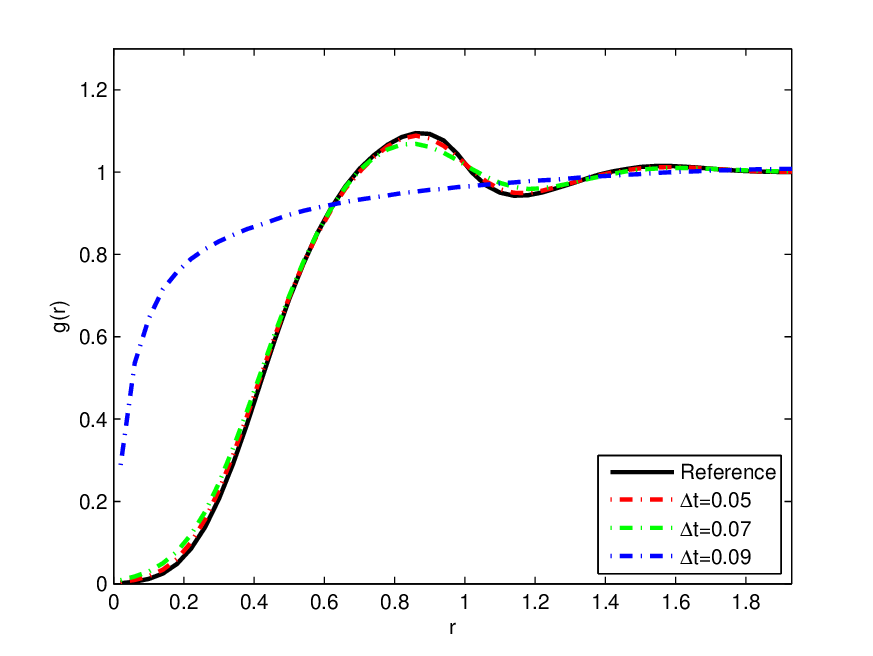}}\hspace{0.2cm}
		\subfloat[M-Shardlow-2]{ \includegraphics[width=0.48\textwidth]{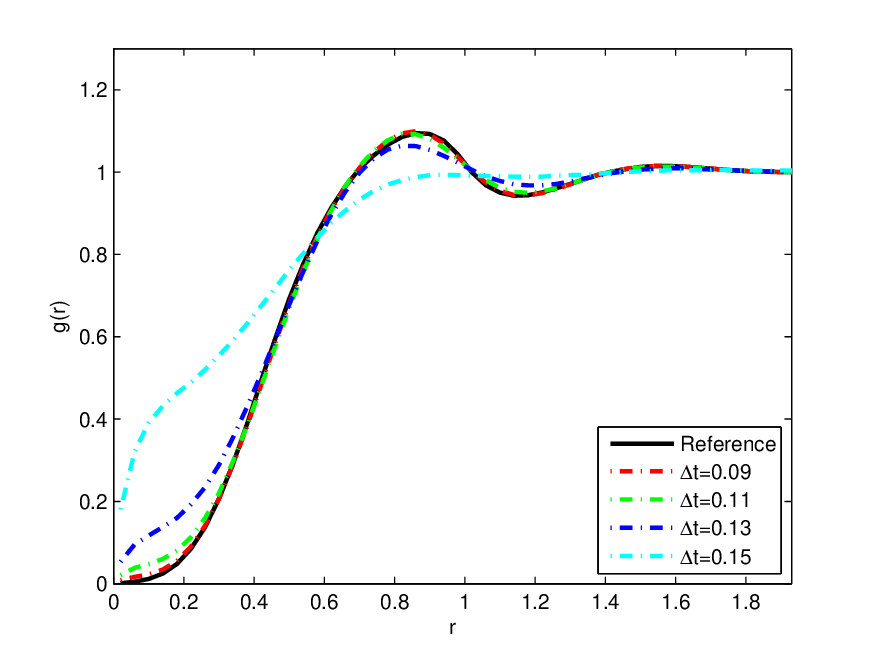}}
	\caption{Comparison of the radial distribution function (RDF), $ g(r) $, for different methods with a large coefficient 
	friction of $ \gamma=40.5 $ and constant repulsion $ a=25 $. In all the methods the stepsize start at $ \Delta t=0.09 $ and RDFs have been plotted until large artifacts observed.
				  				 		The format of the plot is the same as in Figure 3.}\label{fig4}
				  				 \end{figure}

In Figure \ref{fig4} for large friction $ \gamma=40.5 $ and $ a=25 $, the RDFs of the different methods have been depicted. In Shardlow and M-Shardlow-2  methods, for stepsize $ \Delta t=0.11 $ a small deviation from the reference solution can be observed, but in both methods for stepsize $\Delta t=0.13 $ the RDFs are destroyed. For the ABOBA method, plotting the RDFs indicates that  until the stepsize $ \Delta t =0.17 $, a very small artifacts are observed, but around $ \Delta t=0.19 $ the RDFs is destroyed. So for this case, we can see that using a fourth order symplectic method almost did not improve the performance of the Shardlow's method. 
By this figure the high instability of the M-Shadlow-1 method compare to other methods even for stepsize $ \Delta t=0.09 $ is visible.

\begin{figure}[H]
  	\centering
  	\subfloat[Shardlow]{ \includegraphics[width=0.48\textwidth]{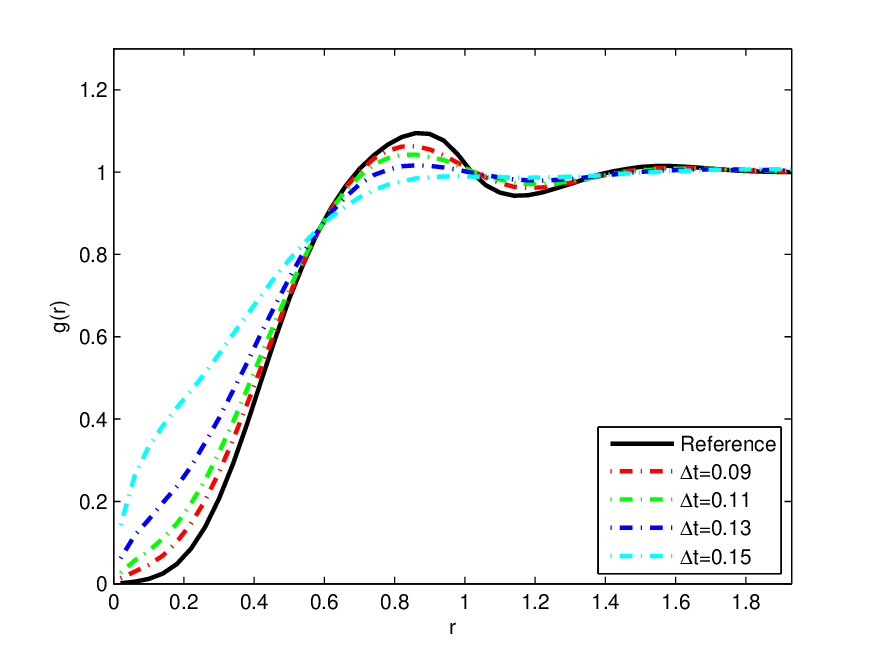}}\hspace{0.2cm}
  	\subfloat[ABOBA]{ \includegraphics[width=0.48\textwidth]{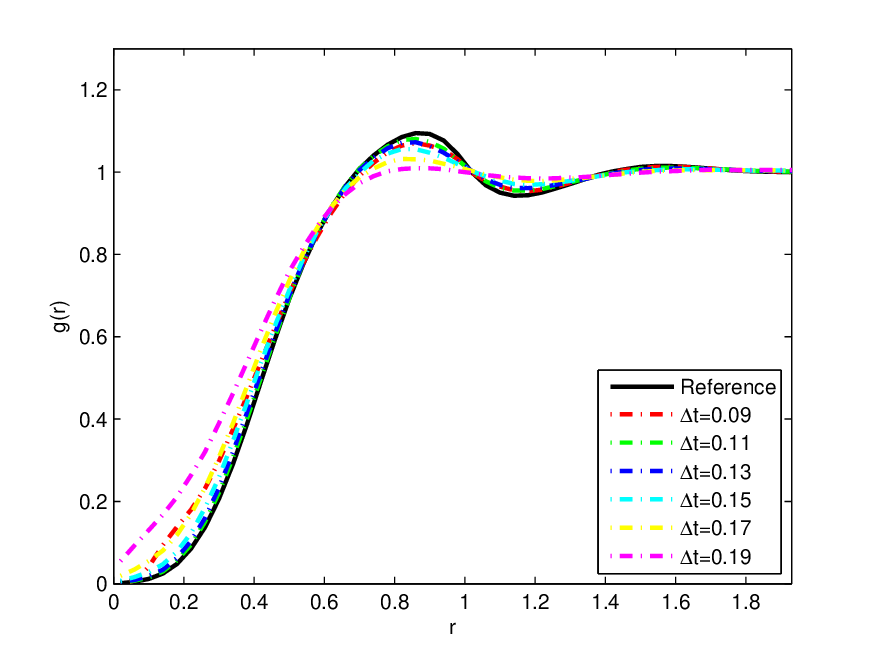}}\\
  	\subfloat[M-Shardlow-1]{ \includegraphics[width=0.48\textwidth]{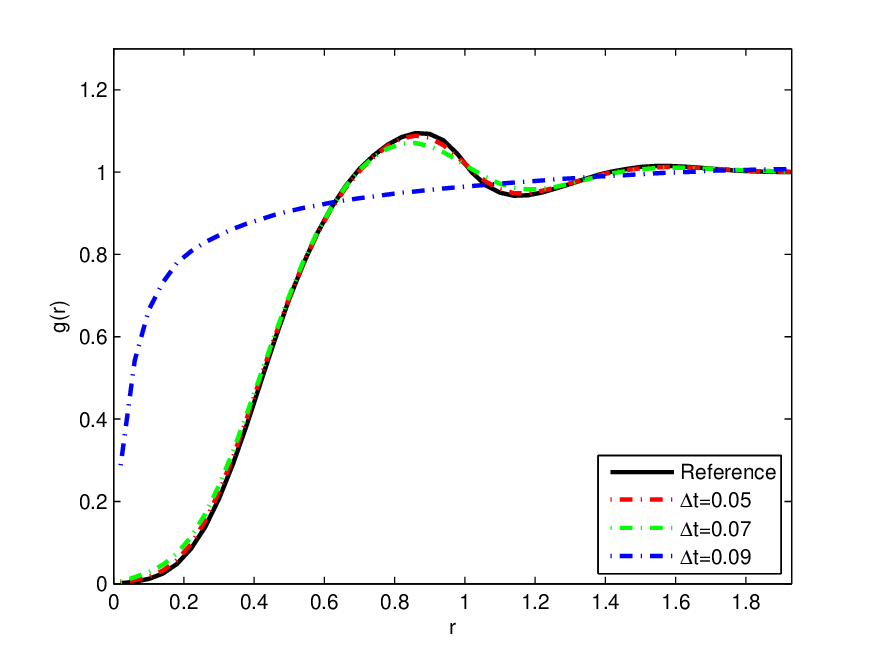}}\hspace{0.2cm}
  	\subfloat[M-Shardlow-2]{ \includegraphics[width=0.48\textwidth]{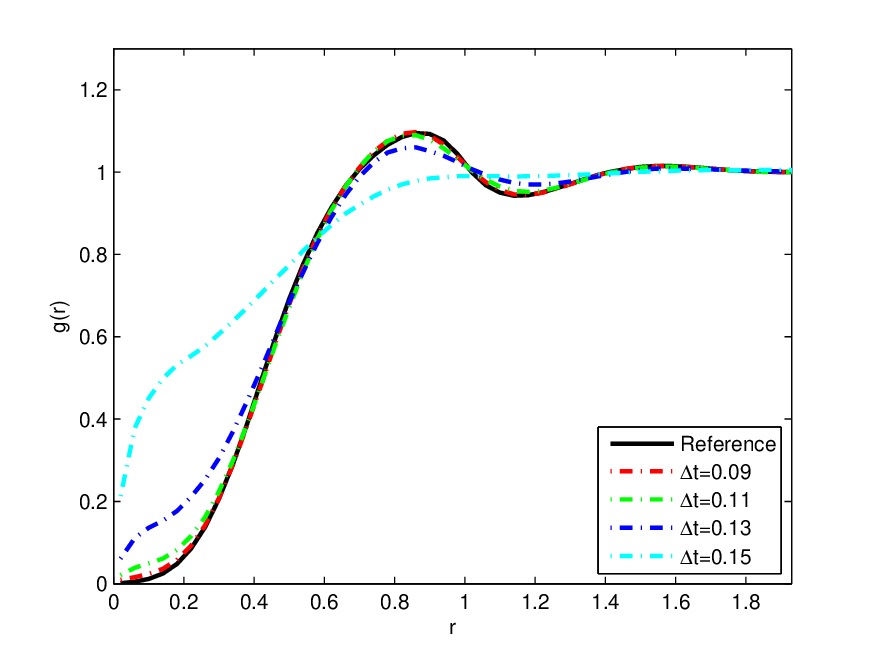}}
  	\caption{Comparison of the radial distribution function (RDF), $ g(r) $, for different methods with a extremely large coefficient friction of $ \gamma=450 $ and $ a=25 $. In all the methods the stepsize start at $ \Delta t=0.09 $ and RDFs have been plotted until large artifacts observed. The format of the plots is the same as in Figure 3.}\label{fig5}
  \end{figure}

The performance for $ (\gamma,a)=(200,25) $ and $ (\gamma,a)=(450,25) $ are very similar. So, we only plotted the RDF for $ (\gamma,a)=(450,25) $  in Figure \ref{fig5}. These results are similar to the results we obtain for $ \gamma=40.5 $, where again the M-Shardlow-2 method show no considerable superiority over Shardlow and ABOBA methods. The RDFs for stepsize smaller than $ \Delta t=0.11 $ for M-Shardlow-2 method is almost indistinguishable from the reference RDF, but at $ \Delta t=0.15 $ it has been completely destroyed. The ABOBA method controls the structure of the system very well for larger increments and even at  $ \Delta t=0.19 $ still ABOBA shows a good performance. To observe the instability of M-Shardlow-1 method, similar to the previous cases, the drawing of the
RDF for this method is also included.  

				  		 \section{Conclusions}	
				  		  In this paper, a modification of Sardlow's method has been presented. This method that we have named M-Shardlow-2 is an improvement of a modified method which was initially proposed by Yamada et.al. \cite{Yamada}, in which the authors replaced the second order velocity Verlet in Sharlow's method with a fourth order symplectic scheme. In this paper, we showed that authors in \cite{Yamada} have used improper parameters to construct their modified method. Then, according to Yoshida's scheme for constructing high order symplectic methods, we propose an improvement of this method. Besides, both numerical and theoretical results verified that the M-Shardlow-2 method demonstrates a third order convergence to the invariant measure.  
				  		  
				  		  On the other hand, although, the M-Shardlow-2 method increase the order of convergence and consequently improve the stability threshold considerably, but testing the numerical efficiency shows that the most disadvantages of the modified M-Shardlow-2 is that the required CPU time for this method that for each step is almost three times more than that of Shardlow and ABOBA methods. Three evaluations of the conservative force in the M-Shardlow-2 method compare to one evaluation of conservative force in Shardlow and ABOBA methods eliminates the benefits of using this method. 
				  		  
				  		  Configurational quantities such as configurational temperature and radial distribution function (RDF) were tested for the M-Shardlow-2 method. The obtained results for RDF in many experiments show that it is not able to outperform the Shardlow and ABOBA methods. Even in some cases, it performs weaker in preserving the structure and dynamical properties of the system compared to Shardlow and ABOBA methods. Overall it seems that using the proposed procedure by the authors in \cite{Yamada} to build high order integrators for DPD systems is useless and has no substantial justifications.   
				
			\end{document}